\documentclass{article}

\RequirePackage{amsmath}

\RequirePackage{hyperref}

\usepackage{arxiv}

\usepackage[utf8]{inputenc} 
\usepackage[T1]{fontenc}    
\usepackage{hyperref}       
\usepackage{url}            
\usepackage{booktabs}       
\usepackage{amsfonts}       
\usepackage{nicefrac}       
\usepackage{microtype}      
\usepackage{lipsum}		

\usepackage{graphicx}
\usepackage{subfigure}
\usepackage[super]{nth}
\usepackage{mathptmx}     
\usepackage{multirow}
\usepackage{algorithm}
\usepackage{algorithmicx}
\usepackage{algcompatible}
\usepackage{multicol}
\usepackage{algpseudocode}
\usepackage[super]{nth}
\usepackage{enumitem}

\newtheorem{lem}{Lemma}

\newtheorem{assm}{Assumption}
\newcommand{\norm}[1]{\Vert #1 \Vert}
\newcommand{\E}{\mathbb{E}}
\newcommand{\Q}{\mathfrak{Q}}
\newcommand{\inner}[2]{\langle#1,#2\rangle}

\newcommand{\exclude}[1]{}
\newcommand{\N}{\mathcal{N}}
\renewcommand{\P}{\mathcal{P}}
\newcommand{\T}{\mathcal{T}}
\newcommand{\Z}{\mathcal{Z}}

\usepackage{tikz}
\usetikzlibrary{arrows, decorations.pathmorphing,decorations.markings,trees,arrows.meta}
\tikzset{
  treenode/.style = {align=center, inner sep=0pt, text centered,
    font=\sffamily},
sqr_x/.style = {treenode, circle, line width=1mm, white, draw=black,
    fill=white, minimum width=1cm, minimum height=1cm},
sqr_xx/.style = {treenode, circle, line width=1.25mm, white, draw=black,
    fill=white, minimum width=0.3cm, minimum height=0.3cm},
  arn_x/.style = {treenode, circle, white, draw=black,
    fill=black, minimum width=0.75cm, minimum height=0.75cm},
  arn_r/.style = {treenode, circle, white, draw=red,
    fill=red, minimum width=0.75cm, minimum height=0.75cm},
  arn_b/.style = {treenode, circle, white, draw=blue,
    fill=blue, minimum width=0.75cm, minimum height=0.75cm},
  arn_xx/.style = {treenode, circle, white, draw=black,
    fill=black, minimum width=0.25cm, minimum height=0.25cm},
  arn_rr/.style = {treenode, circle, white, draw=red,
    fill=red, minimum width=0.25cm, minimum height=0.25cm},
  arn_bb/.style = {treenode, circle, white, draw=blue,
    fill=blue, minimum width=0.25cm, minimum height=0.25cm},
  arn_w/.style = {treenode, circle, white, draw=white,
    fill=white, minimum width=0.125cm, minimum height=0.125cm},
  arn_vsw/.style = {treenode, circle, white, draw=white,
    fill=white, minimum width=0.0001cm, minimum height=0.0001cm},
  arn_label/.style = {treenode, circle, white, draw=black,
    fill=black, minimum width=0.1cm, minimum height=0.1cm},
	CC/.style={-Latex, decorate, decoration={snake}, draw=black, line width=1mm},
	fwd/.style={-Latex, draw=black,line width=1mm},
	infwd/.style={-Latex, draw=black!30!green,line width=1mm},
	fwdds/.style={-Latex, draw=black,line width=1mm, dashed, dash pattern= on 17.5pt},
	Fwdds/.style={-Latex, draw=black,line width=1mm, dashed, dash pattern= on 2pt},
	Fwds/.style={-Latex, draw=black!30!green,line width=1mm,densely dash dot},
	SCC/.style={-Latex, decorate, decoration={snake}, draw=black, line width=1mm, densely dotted},
	arn_c/.style = {treenode, rectangle, white, draw=black,
    fill=white, minimum width=1.9cm, minimum height=0.8cm, line width=1mm, dashed},
    	arn_cs/.style = {treenode, rectangle, white, draw=black,
    fill=white, minimum width=1.6cm, minimum height=0.8cm, line width=1mm, dashed},
        	arn_cvs/.style = {treenode, rectangle, white, draw=black,
    fill=white, minimum width=0.8cm, minimum height=0.8cm, line width=1mm, dashed}
}

\title{Adaptive Partition-based SDDP Algorithms for Multistage Stochastic Linear Programming}

\author{
Murwan Siddig \\
  Department of Industrial Engineering\\
  Clemson University\\
  Clemson, SC, USA 29631  \\
  \texttt{msiddig@clemson.edu} \\
   \And
 Yongjia Song \\
  Department of Industrial Engineering\\
  Clemson University\\
  Clemson, SC, USA 29631  \\
  \texttt{yongjis@clemson.edu} \\
}

\begin{document}
\maketitle

\begin{abstract}
In this paper, we extend the adaptive partition-based approach for solving two-stage stochastic programs with fixed recourse to the multistage stochastic programming setting. The proposed algorithms integrate the adaptive partition-based strategy with a popular approach for solving multistage stochastic programs, the \textit{stochastic dual dynamic programming}, via different tree-traversal strategies in order to enhance its computational efficiency. Our numerical experiments on a hydro-thermal power generation planning problem show the effectiveness of the proposed algorithms.

\end{abstract}


\section{Introduction}
\label{intro}
Multistage stochastic programming is a well-recognized mathematical optimization model for problems that require optimization under uncertainty over time. These problem arise in a variety of applications, such as energy \cite{goel2004stochastic,de2017assessing,homem2011sampling,pereira2005strategic,rebennack2016combining,shapiro2013risk}, finance \cite{ahmed2003multi,ch2007modeling,dupavcova2009portfolio,dupavcova2009asset}, transportation \cite{alonso2000stochastic,fhoula2013stochastic,herer2006multilocation} and sports \cite{pantuso2017football}, among others. In particular, in this paper, we consider a class of multistage stochastic programming models which have applications in long-term hydro-thermal power generation planning, in which one aims to construct an optimal operational strategy under the uncertainty of rainfall volume in order to minimize power generation cost to meet certain demand. The sequential nature in the decision-making structure and the uncertainty in the problem data such as, future (water) inflows, demand, fuel cost, etc., make the hydro-thermal power generation planning problem a classical problem in applications of multistage stochastic programming.

In the stochastic programming literature, there is a good deal of work on how to tackle such problems in a computationally tractable fashion. One typical approach is to approximate the underlying stochastic process governing the uncertainty in the problem data using scenario trees. The result of doing this is two fold:
\begin{enumerate}
\item As the number of decision stages in the planning horizon increases, the increasing scenario-tree size requires an exponential growth of computational resources for solving the corresponding multistage stochastic program~\cite{bellman}.
\item Specialized algorithms that employ decomposition techniques become a necessity to solve the resulting large-scale mathematical programs. 
\end{enumerate}

In this paper we are concerned with the computational efficiency for solving multistage stochastic programming problems based on enhancements to one of the most successful decomposition algorithms for these problems, the \textit{Stochastic Dual Dynamic Prgramming} (SDDP) algorithm \cite{SDDP}. The proposed enhancements are based on the idea of employing adaptive partition-based formulation \cite{song2017,song2015}, which gives a relaxation of the original mathematical program obtained by aggregating variables and constraints according to a partition over the set of scenarios. Using the dual information associated with each scenario, the partition can be refined during the solution process until it yields an optimal solution to the problem. Partition-based strategies have shown to be very effective in the two-stage stochastic programming setting due to the reduced computational effort in generating cutting planes that are adaptively refined in the solution process. The task of integrating adaptive partition-based strategies to the SDDP algorithm poses the following questions:
\begin{enumerate}
 \item How to choose the way for the SDDP algorithm to traverse the scenario tree with aggregated variables and constraints according to the scenario parition in each stage?
 \item For a given tree-traversal strategy, how accurate should the solution be during different phases of the solution process?
 \item How should we refine the scenario partition in each stage?
 \end{enumerate}
Using different strategies to traverse the scenario tree, refining the partition in an \textit{adaptive partition-based SDDP} algorithm can be done in many different ways. We investigate two different strategies, namely \emph{refinement outside of the SDDP algorithm}, and \emph{refinement within the SDDP algorithm}. Moreover, we develop a method which exploits the nature in problems of optimal policies with special structures. This is done by incorporating partition-based strategies only to a selected subset of stages in the planning horizon, while the standard approach in SDDP is applied to other stages.

The rest of this paper is organized as follows: in Section \ref{sec:preliminaries} we introduce a general multistage stochastic programming formulation and provide an overview of the theoretical background in adaptive partition-based strategies and the SDDP algorithm. In Section \ref{sec:PartitionSDDP}, we demonstrate how the adaptive partition-based approach can be used in the multistage setting and provide the ingredients of our proposed algorithms. In Section \ref{sec:ImplmntDetails} we describe in details three types of adaptive partition-based SDDP algorithms: \textit{Refinement outside SDDP},  \textit{Refinement within SDDP} and \textit{Adaptive partition-based SDDP with structured cut generation policies}. In Section \ref{sec:results} an extensive computational analysis is presented, and the proposed approach is compared with alternative approaches in terms of their computational performance. Finally, in Section \ref{sec:conclusions} we conclude with some final remarks.  
\section{Preliminaries on Multistage Stochastic Linear Programs and Decomposition Schemes}
\label{sec:preliminaries}
In stochastic programming (SP) models, the underlying data uncertainty involved in the problem is characterized by a random vector $\xi$ with known probability distribution. In its most simple form, two-stage SP, two kinds of decisions are involved: \textit{first-stage} decisions $x_1$ that are made prior to the realization of random vector $\xi$, and \textit{second-stage} recourse decisions $x_2 := x_2(\xi)$ that are made after observing the realization of random vector $\xi$. A two-stage stochastic linear program (2SLP) can be formulated as follows~\eqref{eq:2stageSP}:
\begin{equation}
\label{eq:2stageSP}
\min_{\underset{x_1\in \mathbb{R}^{n_1}_+}{A_1x_1=b_1}} c^\top_1  x_1+
\E_\xi\left[
\min_{\underset{x_2\in \mathbb{R}^{n_2}_+}{B_2 x_1 +
A_2 x_2=b_2}}c_2^\top x_2 \right] \,,
\end{equation}
where $\xi = (c_2, B_2, A_2, b_2)$, and the expectation $\E_\xi[\cdot]$ is taken with respect to the probability measure of random vector $\xi \in \Xi$. 

Multistage stochastic linear programs (MSLPs) provide an explicit framework which generalizes the 2SLP for multiple stages of sequential decision making under uncertainty. In a planning horizon of $T$ stages, the dynamic realization of uncertainty is typically modeled as a stochastic process $(\xi_1, \xi_2, \dots \xi_T)$ such that, $\xi_1$ is deterministic, and for each $t = 2,3,\ldots, T$, $\xi_t \in \Xi_t$ is random vector that will be realized in stage $t$. The history of this stochastic process up to time $t$ is denoted by $\xi_{[t]} := (\xi_1, \dots \xi_t)$. The nested form of an MSLP can be expressed as:
\begin{equation}
\label{eq:mlsp}
\min_{\underset{x_1\in \mathbb{R}^{n_1}_+}{A_1x_1=b_1}} c_1^\top  x_1+
\E_{|\xi_{[1]}}\left[
\min_{\underset{x_2\in \mathbb{R}^{n_2}_+}{B_2x_1 +
A_2x_2=b_2}}c_2^\top x_2 + \E_{|\xi_{[2]}}\left[\cdots+\E_{|\xi_{[T-1]}}\left[\min_{\underset{x_T\in \mathbb{R}^{n_T}_+}{B_Tx_{T-1} + A_Tx_T=b_T}}c_T^\top  x_T
\right]\right]\right] \,,
\end{equation}
where some (or all) data $\xi_t=(c_t,B_t,A_t,b_t)$ can be subject to uncertainty for $t=2,\ldots,T$. The expectation $ \E_{|\xi_{[t]}}[\cdot]$ is taken with respect to the conditional probability measure of the random vector $\xi_{t+1}$. The sequence of decisions ($\bar{x}_1, \bar{x}_2, \dots , \bar{x}_t$), denoted by $\bar{x}_t= \bar{x}_t(\xi_{[t]}), \forall \; t=1, 2, \dots, T$ is referred to as an implementable policy for problem~\eqref{eq:mlsp}. Such policy provides a decision rule at every stage $t$ based on the realization of the data process up to time $t$. The aim of an MSLP is to find an optimal policy to~\eqref{eq:mlsp}. 

The question of how to construct scenarios to induce a decision policy and measure its quality is beyond the scope of this paper. In this paper, we assume that a scenario tree $\T$ is given, where a finite number of realizations is available for each $\xi_t, \ \forall \; t = 2,3,\ldots, T$. As such,~\eqref{eq:2stageSP} and~\eqref{eq:mlsp} can be written as large-scale linear programs, known as the \textit{deterministic equivalent programs} (DEP). This is, for instance, the case in which~\eqref{eq:mlsp} is a sample average approximation (SAA) of an original MSLP where the underlying random vector $\xi_t$ in each stage $t$ follows a continuous probability distribution. We refer the reader to \cite{Shapiro_2011} for a discussion on the relationship between an SAA problem and the original MSLP with a continuous distribution.

\subsection{Preliminaries on Two-stage Stochastic Linear Programs}
Consider solving the following DEP for the 2SLP~\eqref{eq:2stageSP} with a set of scenarios indexed by $ N = \{1, 2, \dots |\Xi|\}$:
\begin{equation}
\label{eq:master-problem}
\begin{array}{llll}
\displaystyle z^* = & \underset{x_1 \in \mathbb{R}_{+}^{n_1} }{\min}& \; c^\top_1 x_1+ \Q(x_1)\\ 
& \mbox{s.t.} & A_1x_1=b_1
\end{array}
\end{equation}
where the second-stage cost $\Q(x_1) := \E_\xi[Q(x_1,\xi)]$, such that
\begin{equation}
\label{eq:subproblem}
Q(x_1,\xi) = \min_{x_2 \in \mathbb{R}_{+}^{n_2}} \left\{c_2^\top x_2 \; | \; B_2 x_1 + A_2 x_2=b_2 \right\}.
\end{equation}
For simplicity, we assume that each scenario happens with an equal probability $1/|\Xi|$. It can be easily seen that, as the number of scenarios grows, the DEP becomes computationally challenging to solve. Nevertheless, the DEP has a special structure that lends itself to decomposition techniques developed to solve large-scale LPs. These include variants of \textit{Benders decomposition} (also called the L-shaped method \cite{Lshaped}) for the two-stage setting, which generalizes to \textit{nested Benders decomposition} for MSLPs.

\subsubsection{Benders decomposition and inexact oracles}
\label{subsec:Benders}

Standard Benders decomposition consists of iteratively solving an approximate problem referred to as the \textit{master problem}, where the second-stage cost function $\Q(\cdot)$ in~\eqref{eq:master-problem} is approximated by a cutting-plane lower approximation, and a subproblem for each scenario $k \in N$. The master problem provides a candidate solution $x_1$, which is an optimal solution to the master problem under the current approximation of $\Q(\cdot)$, denoted by $\check \Q(\cdot)$. This solution $x_1$ will then be evaluated by solving the subproblems. The cutting plane approximation for $\Q(\cdot)$ is improved by providing \textit{feasibility cuts} and/or \textit{optimality cuts} generated using information obtained from solving the subproblems. We refer the reader to \cite{birge2011} and \cite{kall1994} for a detailed discussion on the topic.

From an abstract viewpoint coming from nonsmooth optimization, such decomposition schemes can be seen as variants of Kelley's cutting plane method for solving convex programs \cite{kelley1960}. In the terminology of convex programming, an \textit{oracle} is referred to as a routine that returns the function value information and the subgradient information at a given point $x_1$. If the information provided by the routine is accurate, we call the oracle an \textit{exact oracle}. Otherwise, the oracle is called \textit{inexact oracle}. In the context of the 2SLP, given by the DEP~\eqref{eq:master-problem}, oracles constructed by solving the scenario subproblem~\eqref{eq:subproblem} for each $k \in N$ correspond to exact oracles, yielding exact function value $\Q(x_1)$ and subgradient information from the optimal dual solutions.

It is intuitively clear that when the number of scenarios $|N|$ is large, the exact oracle could be computationally expensive. On the other hand, it may be beneficial to use a relatively coarse oracle at the beginning of the algorithm, which is just used to get the process "warm-started"; the exact oracle is used afterwards whenever it is necessary. Oliveira and Sagastizabal \cite{de2014level} formalize this idea, and introduce the concept of inexact oracle with on-demand accuracy for a generic (nonsmooth) convex optimization problem. Inexact oracles can be constructed in many different ways in the context of 2SLPs, van Ackooij, Oliveira, and Song \cite{song2017} study an inexact oracle defined by scenario partitions, which we explain in detail next. To that end, for reasons that will become clear soon, throughout this paper we make the following assumption:
\begin{assm}
\label{assm:fixed-recourse-cost}
\underline{Fixed recourse and fixed cost vector}. We assume that the recourse matrix $A_t$ and the cost vector $c_t$ for $t=2, \dots, T$, are the same for all realizations of $\xi_t$. Hence, the uncertainty in each stage is characterized only by $\{(B_{t,k},b_{t,k})\}_{k=1}^{|\Xi_t|}$.
\end{assm}
As a result of Assumption~\ref{assm:fixed-recourse-cost}, the dual feasible region for the second-stage problem~\eqref{eq:subproblem} is the same for all scenarios: $\Pi := \{\pi \; | \; A^\top_2 \pi \leq c_2 \}$.

\subsubsection{Partition-based inexact oracles} 
\label{subsec:2stage-partition}
A \textit{partition} $\N = \{\P_1, \P_2, \dots, \P_L\}$ of the scenario set $N$ is a collection of nonempty subsets of scenarios such that, $\P_1 \cup \P_2 \cup \dots \P_L = N$, and $\P_i \cap \P_j = \emptyset , \; \forall \; i, j \in \{1, 2, \dots, L\}, i \neq j$. Each of these subsets is called a scenario cluster. Given a partition $\N$, the second-stage cost function can be alternatively written as:
\begin{equation*}
\Q(x_1) = \frac{1}{|\Xi|} \sum_{k \in N} Q(x_1,\xi_k) = \frac{1}{|\Xi|} \sum_{j=1}^{L}\sum_{k \in \P_j} Q(x_1,\xi_k).
\end{equation*}
It is clear that aggregating the scenarios within each cluster defines a relaxation for~\eqref{eq:subproblem}. That is, for each $j=1, \dots L$, define ${\bar{b}}^{\P_j}_2 := \sum_{k \in \P_j} b_{2}^k$, and ${\bar{B}}^{\P_j}_2  := \sum_{k \in \P_j} B_{2}^k$, then $\forall \; x_1 \in X_1$ we have that $\sum_{k \in \P_j} Q(x_1,\xi_k) \geq Q^{\P_j} (x_1)$ where 
\begin{equation}
\label{eq:2stage-partition}
Q^{\P_j} (x_1) := \min_{x_2 \in \mathbb{R}_{+}^{n_2}} \left\{c^\top_2 x_2 \mid  {\bar{B}}^{\P_j}_2  x_1+ A_2 x_2 = {\bar{b}}^{\P_j}_2 \right\}.
\end{equation}
Hence, $\Q(\cdot)$ can be lower approximated by:
\begin{equation}
\label{eq:partition-2cost}
\Q(x_1) = \frac{1}{|\Xi|} \sum_{j=1}^{L}\left[\sum_{k \in \P_j} Q_k(x_1) \right] \geq  \Q^\N (x_1) := \frac{1}{|\Xi|} \sum_{j=1}^{L} Q^{\P_j} (x_1).
\end{equation}
Note that, problem~\eqref{eq:subproblem} is referred to as the \textit{scenario-based subproblem}, and problem~\eqref{eq:2stage-partition} is referred to as the \textit{partition-based subproblem}. 
\paragraph{Partition refinement.} The inexactness of a partition-based oracle is given by the gap between $\Q(x_1)$ and $\Q^\N(x_1)$, which can be reduced by \textit{partition refinement}. Specifically, we say that $\N'$ is a refinement of $\N$, if $|\N'| > |\N|$, and $\forall \; \P' \in \N', \; \exists \; \P \in \N$ such that $\P' \subseteq \P$, i.e., $\P'$ is obtained by subdividing some cluster(s) of $\N$. It is clear that after the partition refinement, $\Q(x_1) \geq \Q^{\N'}(x_1) \geq \Q^{\N}(x_1)$. Song and Luedtke \cite{song2015} propose a refinement strategy driven by the optimal dual solutions to~\eqref{eq:subproblem}, which is motivated by the following result:
\begin{lem}
Let $\N = \{\P_1, \P_2 , \dots \P_L\}$ be a \textit{partition}, and $x_1 \in X_1$. Then  $\Q(x_1) = \Q^\N(x_1)$ if for each $j =1, 2, \dots L$ there exists dual optimal solutions $\pi_j$ for problem~\eqref{eq:2stage-partition} such that $\pi_j \in \Pi^*_k(x_1), \; \forall \; k \in \P_j$, where $\Pi^*_k(x_1)$ is the set of dual optimal solutions to the subproblem~\eqref{eq:subproblem} at a given $x_1$. 
\end{lem}
Therefore, given a solution $x_1 \in X_1$, if $\Q(x_1) > \Q^{\N}(x_1)$, one could refine partition $\N$ into $\N'$ by subdividing each $\P_j, \; \forall j = 1, 2, \dots, L$ according to the optimal dual solutions $\pi_k$ for each scenario $k \in \P_j$ such that after the partition is refined, $\Q^{\N'}(\cdot)$ locally matches with $\Q(\cdot)$ at a given solution $x_1$. We refer the reader to \cite{song2017} and \cite{song2015} for a detailed discussion on the topic.

\subsection{Preliminaries on Multi-stage Stochastic Linear Programs}
We now turn our attention to decomposition schemes used in solving MSLPs. The decomposition schemes for the 2SLP~\eqref{eq:2stageSP} are somewhat readily accessible due to the static nature of the decision structure in the problem. Decisions in the two-stage setting can be viewed as static in the sense that, a supposedly optimal decision was made at a previous point in time and recourse actions are taken after all of the uncertainty is resolved. This is not necessarily the case in the MSLP~\eqref{eq:mlsp} -- now that decisions are made sequentially over time, based on the information available at every stage. When $T > 2$, the expectation taken in~\eqref{eq:mlsp} at time $t=T-1$ is essentially the same as the second-stage cost in~\eqref{eq:2stageSP}, except that in~\eqref{eq:mlsp} the expectation is taken with respect to the conditional probability distribution based on the history of the stochastic process up to time $t=T-1$, i.e., $\xi_{[T-1]}$. Therefore, if one were to define a function at every stage $t =T-1, T-2, \dots 1$, in order to capture the future cost in an exact way, i.e., the same way that $\Q(x_1)$ in~\eqref{eq:master-problem} captures the second-stage cost, we need to calculate the value of those functions recursively going backward in time. In this case, the value of such functions will depend on both $x_t$ and the multivariate variables $\xi_{[t]}$. Hence, depending on the size of the scenario tree, this dependency structure could make the computation of the value of those functions very difficult or even impossible. This motivates our next assumption. We refer the reader to \cite{shapiro2009lectures} for a further discussion on the topic. 
\begin{assm}
\label{assm:stagewise-indp}
\underline{Stage-wise independence}. We assume that the stochastic process $\{\xi_t\}$ is stage-wise independent, i.e., $\xi_t$ is independent of the history of the stochastic process up to time $t-1$, for $t=1, 2, \dots T$, which is given by $\xi_{[t-1]}$.
\end{assm}
Hence, the conditional expectation $\E_{|\xi_{[t-1]}}[\cdot]$ can simply be written as the (unconditional) expectation $\E[\cdot]$ with respect to $\xi_t$. Moreover, note that under the \textit{stage-wise independence} assumption, the optimal policy is a function of $x_t$ alone. This allows for problem~\eqref{eq:mlsp} to be formulated using the following dynamic programming equations:
\begin{equation}
\label{eq:cost2go}
Q_t(x_{t-1},\xi_{t}):= \left\{
\begin{array}{llll}
\displaystyle\min_{x_t \in \mathbb{R}^{n_t}_+} & c_t^\top  x_t + \Q_{t+1}(x_t)\\
\mbox{s.t.} &A_tx_t=b_t-B_tx_{t-1} \,,
\end{array}
\right.
\end{equation}
where $\Q_{t+1}(x_t):= \E[Q_{t+1}(x_t,\xi_{t+1})]$ for $t=T-1,\ldots 1$,
and $\Q_{T+1}(x_{T}):=0$. The first-stage problem becomes
\begin{equation}
\label{eq:1stagepbm}
\left\{
\begin{array}{llll}
\displaystyle \min_{x_1\in \mathbb{R}^{n_1}_+}& c_1^\top x_1 + \Q_2(x_1)\\
\mbox{s.t.}& A_1x_1=b_1 \,.
           \end{array}
           \right.
\end{equation}
$\Q_{t+1}(\cdot)$ in~\eqref{eq:cost2go} is referred to as the \textit{cost-to-go} function. Note that, if the number of scenarios per stage is finite, the cost-to-go functions are convex piecewise linear functions \cite[Chap. 3]{Ruszczynski_Shapiro_2009b}. Before we discuss the solution approaches developed to solve problem~\eqref{eq:1stagepbm}, it is important to address boundedness and feasibility assumptions.

\paragraph{Boundedness:} The cost-to-go functions $\Q_{t+1}(\cdot)$'s are defined as the optimal future cost at stage $t$, for $t=1, 2 \dots, T-1$. It may happen that for some feasible $x_t \in X_t$ and a scenario $\xi_{t,k} \in \Xi_t$, the \textit{stage}-$(t+1)$ subproblem is unbounded from below, i.e., $Q_{t+1}(x_{t}, \xi_{t+1,k}) = -\infty$. This is a somewhat pathological and unrealistic situation meaning that for such feasible $x_t$, there exist a positive probability, by which one can reduce the future cost indefinitely. One should make sure at the modeling phase that this does not happen.
\begin{assm}
\label{assm:boundedness}
\underline{Boundedness}. We assume that the cost-to-go function $\Q_{t+1}(\cdot)$ at every stage is finite valued. i.e., $Q_{t+1}(x_t, \xi_{t+1,k}) > -\infty , \ \forall \xi_{t,k} \in \Xi_t, \ \forall t=1, 2 \dots, T-1$.
\end{assm}
\paragraph{Feasibility:} If for some $x_t \in X_t$ and scenario $\xi_{t,k} \in \Xi_t$ problem~\eqref{eq:cost2go} is infeasible, the standard practice is to set $Q_{t+1}(x_t, \xi_{t,k}) = +\infty$ such that $x_t$ cannot be an optimal solution of the \textit{stage}-$(t+1)$ subproblem. It is said that the problem has relatively complete recourse if such infeasibility does not happen.
\begin{assm}
\label{assm:relatively-complete-recourse} 
\underline{Relatively complete recourse}. We assume that $\forall \; x_t \in X_t$ and $\xi_{t,k} \in \Xi_t$, $\exists \; x_{t+1} \in \mathbb{R}^{n_{t+1}}_{+}$ such that, problem~\eqref{eq:cost2go} is feasible. 
\end{assm}

\subsubsection{Stochastic Dual Dynamic Programming (SDDP)} 
\label{subsec:sddp}
As one of the most popular algorithms for solving MSLPs, the SDDP algorithm \cite{SDDP} draws influence from the backward recursion techniques developed in dynamic programming \cite{bellman}. More importantly, under Assumption \ref{assm:stagewise-indp}, it provides an implementable policy not only for the approximation problem, but also for the true problem through decision rules induced by the approximate cost-to-go functions.

Leveraging the dynamic equations~\eqref{eq:cost2go}, the SDDP algorithm alternates between two main steps: a forward simulation (\textit{forward pass}) which evaluates the current policy obtained by the current approximation for the cost-to-go functions $\check \Q_{t+1}(\cdot)$, and a backward recursion (\textit{backward pass}) to improve $\check \Q_{t+1}(\cdot)$, for $t=2, \dots T$. After the forward step, a statistical upper bound $\overline{z}$  for the optimal value is determined, and after the backward step a lower bound $\underline z$ is obtained. We summarize the two main steps below and refer the reader to \cite{SDDP} for a more detailed discussion on the topic.\\

\noindent\underline{Forward step.} Consider taking a sample of $\mathcal{J}$ scenarios of the stochastic process, which is denoted by $\{\xi^j\}_{j\in \mathcal{J}}$, with $|\mathcal{J}| \ll |\Xi_1| \times |\Xi_2| \times \dots \times |\Xi_T|$ and $\xi^j = (\xi_1^j, \ldots,\xi_T^j)$. Let $\check \Q_{t+1}(\cdot)$ be the current approximation of the cost-to-go function $ \Q_{t+1}(\cdot)$ at stage $t$. In order to evaluate the quality of the decision policy induced by $\check \Q_{t+1}(\cdot)$ for $t=1, 2, \dots, T$, trial decisions $\bar x_t=\bar x_t(\xi_{[t]})$, $t=1,\ldots,T$, are computed recursively going forward with $\bar x_1$ being an optimal solution of~\eqref{eq:1stagepbm} with $\check \Q_{2}(\cdot)$, and $\bar x_t$ being an optimal solution of
\begin{equation}
\label{eq:forwt}
\underline{Q}_{t}(\bar{x}_{t-1},\xi_{t}^j):= \left\{
\begin{array}{llll}
\displaystyle \min_{x_{t}\in \mathbb{R}^{n_t}_+}& c_{t}^\top  x_{t} + \check \Q_{t+1}(x_{t})\\
\mbox{s.t.} & A_{t}x_{t}=b_{t} - B_{t}\bar x_{t-1}
\end{array}
\right.
\end{equation}
where $\check \Q_{T+1}(\cdot):= 0$ and $\underline{Q}_{t}(\cdot)$ is an lower approximation for $Q_t(\cdot)$. The value
\begin{equation}
\label{zsup}
z(\xi^j) = \sum_{t=1}^T c_t^\top \bar x_t(\xi_{[t]}^j), \; \forall \; j \in \mathcal{J},
\end{equation}
as well as 
\begin{equation}
\label{eq:forwt}
\begin{array}{llll}
\tilde  z= \frac{1}{|\mathcal{J}|}\sum_{j \in \mathcal{J}}z(\xi^j), \; \forall \; j \in \mathcal{J}, \\
\sigma^2=\frac{1}{|\mathcal{J}|}\sum_{j \in \mathcal{J}}(z(\xi^j)-\tilde z)^2, \; \forall \; j \in \mathcal{J},
\end{array}
\end{equation}
are computed, with $\tilde{z}$ being the sample average and $\sigma^2$ being the sample variance of $z(\cdot)$. Note that $\bar x_t(\xi_{[t]}^j)$ is a feasible and implementable policy for problem~\eqref{eq:mlsp}.

The sample average $\tilde{z}$ provides an unbiased estimator for an upper bound of the optimal
value of~\eqref{eq:mlsp}, which is given by
\begin{equation}
\label{zsup}
\overline{z}=\E\left[\sum_{t=1}^T c_t^\top \bar
x_t(\xi_{[t]}).\right]
\end{equation}
Additionally, $\tilde z +1.96\tilde \sigma/\sqrt{|\mathcal{J}|}$ provides a statistical upper bound for the optimal value of~\eqref{eq:mlsp} with $95\%$ confidence level. These bounds can be used in a possible stopping criterion whenever $\tilde z + 1.96\tilde \sigma/\sqrt{|\mathcal{J}|} -\underline{z}\leq \epsilon$, for a given tolerance $\epsilon>0$. We refer to \cite[Sec.3]{Shapiro_2011} for a discussion on this subject. 
\noindent \underline{Backward step.} 
Given the trial decisions $\bar x_t=\bar x_t(\xi_{[t]})$ obtained in the \textit{forward step} and an approximation of the cost-to-go function $\check{\Q}_{t+1}(\cdot)$, for $t =1, 2,\ldots, T$; exploiting the fact that $\Q_{T+1}(\cdot) := 0$, at stage $t=T$ the following problem is solved for each $\xi_T  \in \{(B_{T,k},b_{T,k})\}_{k=1}^{|\Xi_T|}$:
\begin{equation}
\label{eq:cpT}
\underline{Q}_T(\bar x_{T-1},\xi_T)= \left\{
\begin{array}{llll}
\displaystyle\min_{x_T\in \mathbb{R}^{n_T}_+} & c_T^\top  x_T\\
\mbox{s.t.} & A_Tx_T=b_T - B_T\bar x_{T-1}.
\end{array}
\right.
\end{equation}
Let $\bar \pi_T= \bar \pi_T(\xi_T)$ be an optimal dual solution of problem~\eqref{eq:cpT}. Then define $\alpha_T:= \E[b_T^\top \bar \pi_T]$ and $ \beta_T:=
-\E[B_T^\top \bar \pi_T]\in \partial \Q_T(\bar x_{T-1})$ such that 
\[
q_T(x_{T-1}):= \beta_T^{\top} x_{T-1} + \alpha_T \quad = \quad
\Q_T(\bar x_{T-1})+ \inner{\beta_T}{x_{T-1} - \bar x_{T-1}}\,,
\] is a lower cutting-plane approximation for $\Q_{T}(x_{T-1})$ satisfying
\[ 
\Q_{T}(x_{T-1}) \geq q_T (x_{T-1}) \quad \forall \; x_{T-1}\,,
\]
with $q_T(\cdot)$ being a supporting hyperplane for $\Q_{T}(\cdot)$, i.e., $\Q_{T}(\bar x_{T-1}) = q_T(\bar x_{T-1})$. 
This linear cutting-plane approximation is added to the collection of supporting hyperplanes of $\Q_T(\cdot)$ by 
replacing $\check \Q_T(\cdot)$ with $\check\Q_T(x_{T-1}):= \max\{\check \Q_T(x_{T-1}),q_T(x_{T-1})\}
$. That is, the cutting-plane approximation for $\check \Q_T(\cdot)$ is constructed from the maximum of a collection $J_T$ of cutting-plane approximation:
\[
\check\Q_T(x_{T-1}) = \max_{j \in J_T}\left\{\beta_{T}^{j\top} x_{T-1} +
\alpha_{T}^j \right\}\,.
\]
For $t = T-1,T-2,\ldots, 2$, we update $\check \Q_t(\cdot)$ in the same spirit, the following problems are solved 
\begin{equation}\label{cpt-1}
\underline{Q}_{t}(\bar x_{t-1},\xi_{t})=\left\{
\begin{array}{llll}
\displaystyle \min_{x_{t}\in \mathbb{R}^{n_t}_+}& c_{t}^\top  x_{t} + \check \Q_{t+1}(x_{t})\\
\mbox{s.t.}& A_{t}x_{t}=b_{t} - B_{t}\bar x_{t-1}
\end{array}
\hspace{-0.1cm}
\right. \; \equiv \; \left\{
\begin{array}{llll}
\displaystyle \min_{(x_{t},r_{t+1})\in \mathbb{R}^{n_t}_+\times \mathbb{R}}& c_{t}^\top  x_{t} + r_{t+1}\\
\mbox{s.t.}&  A_{t}x_{t}=b_{t} - B_{t}\bar x_{t-1}\\
 & \beta_{t+1}^{j\top} x_{t} + \alpha_{t+1}^j \leq r_{t+1},\; \; j \in J_{t+1}.
\end{array}
\right.
\end{equation}
$\forall \; \xi_t  \in \{(B_{t,k},b_{t,k})\}_{k=1}^{|\Xi_t|}$, $\forall \; t=T-1, \dots, 2$. Let $\bar \pi_t=\bar \pi_t(\xi_t)$ and $\bar \rho_t^{j} = \bar \rho_t^j(\xi_t)$ be optimal dual multipliers associated with constraints $A_{t}x_{t}=b_{t} - B_{t}\bar x_{t-1}$ and $\beta_{t+1}^{j\top} x_{t} + \alpha_{t+1}^j \leq r_{t+1}$, respectively. Then the linear cutting-plane approximation 
\[
q_t(x_{t-1}):= \beta_t^{\top} x_{t-1} + \alpha_t \quad = \quad
\E[\underline{Q}_{t}(\bar x_{t-1},\xi_{t})] +
\inner{\beta_t}{x_{t-1} - \bar x_{t-1}}
\]
of $\Q_t(\cdot)$ is constructed with
 \begin{equation}\label{cutcoeff}
 \alpha_t:= \E\left[b_t^\top \bar \pi_t + \sum_{j\in J_{t+1}} \alpha_{t+1}^j \bar \rho_t^j\right] \quad\mbox{and}\quad
 \beta_t:= -\E[B_t^\top \bar \pi_t] \in \E[\partial\underline{Q}_{t}(\bar x_{t-1},\xi_{t})]
 \end{equation} such that
$ \Q_{t}(x_{t-1}) \geq q_t (x_{t-1})\;\; \forall \; x_{t-1}$. 
Note that the above inequality holds when $t=T-1$ since $\check \Q_T(\cdot)$ approximates $\Q_T(\cdot)$ from below; then $\underline{Q}_{T-1}(\bar x_{T-2},\xi_{T-1})\leq {Q}_{T-1}(\bar x_{T-2},\xi_{T-1})$ implying that $q_{T-1}(\cdot)$ underestimates $\Q_{T-1}(\cdot)$. The result for stage $t$ follows by using a backward induction argument from $T, T-1, \dots, t$.

Once the cutting-plane $q_{T-1}(x_{t-1})$ is computed, then the current approximation $\check\Q_t(x_{t-1})$ is updated at stage $t$ by: $\check\Q_t(x_{t-1})=\max\{\check\Q_t(x_{t-1}),q_t(x_{t-1})\}$\,. It is worth noting that, unlike in stage $T$ where the value of the cost-to-go function is precisely known ($\Q_{T+1} = 0$), the linear approximation given (in early iterations) by $q_{t}(\cdot), \; \forall t=T-1, \dots, 1$, might be a strict under-estimator of $\Q_t(\cdot)$ for all feasible $x_{t-1}$. In other words, $q_t(\cdot)$ might only be a cutting plane but not necessarily a supporting hyperplane. 

Finally, at $t=1$, the following LP is solved
\begin{equation}\label{cp1}
\underline{z}=\left\{
\begin{array}{llll}
\displaystyle \min_{x_{1}\in \mathbb{R}^{n_1}_+}& c_{1}^\top  x_{1} + \check \Q_{2}(x_{1})\\
\mbox{s.t.}& A_{1}x_{1}=b_{1}
\end{array}
\right. \quad \equiv \quad \left\{
\begin{array}{llll}
\displaystyle \min_{(x_{1},r_{2})\in \mathbb{R}^{n_1}_+\times \mathbb{R}}& c_{1}^\top  x_{1} + r_{2}\\
\mbox{s.t.}&  A_{1}x_{1}=b_{1}\\
 & \beta_{2}^{j\top} x_{1} + \alpha_2^j \leq r_{2},\quad j\in J_2\,.
\end{array}
\right.
\end{equation}
The value $\underline{z}$ provides a lower bound for the optimal value of~\eqref{eq:mlsp}. The updated $\check \Q_t(\cdot)$ for $t=2,\ldots, T$, can be used to induce an implementable policy. The convergence analysis of the method can be found in \cite{chen1999convergent,linowsky2005convergence,philpott2008convergence,shapiro2011analysis}.

In the literature, extensions to SDDP and several other algorithms have been proposed. Guigues \cite{guigues2017inexact,guigues2018inexact} proposes an inexact variant of the SDDP algorithm where some or all problems to be solved in the forward and backward passes are solved approximately. In \cite{asamov2018regularized,sen2014multistage,van2019level} the authors deploy regularization techniques that stabilizes iterates during the forward pass of the algorithm. More specifically, in \cite{asamov2018regularized,sen2014multistage} regularization is performed by adding a quadratic term to the objective function of the subproblems, whereas \cite{van2019level} define trial points either as normal solutions of the LP subproblems solved by SDDP or as specific points in the level sets of cutting planes approximating for the cost-to-go functions. In \cite{leclere2018exact} the authors present a dual SDDP algorithm that yields a converging exact upper bound for the optimal value, \cite{baucke2017deterministic} studies a deterministic sampling algorithm and \cite{homem2011sampling} studies sampling strategies and stopping criteria for the SDDP algorithm. 
\section{Adaptive Partition-based SDDP for Multistage Stochastic Linear Programs}
\label{sec:PartitionSDDP}
In this section we demonstrate how the adaptive partition-based strategies discussed in Section~\ref{subsec:2stage-partition} can be used to approximate the cost-to-go functions $\Q_{t+1}(\cdot)$ in~\eqref{eq:cost2go}, and present the ingredients of our proposed algorithms which integrate the adaptive partition-based strategies \cite{song2017,song2015} and the SDDP algorithm \cite{SDDP}.
\subsection{Adaptive Partition in the Backward Step}
\label{subsec:backward-partition}
Let $\bar{x}_t = \bar{x}_t (\xi_t), t=1,2,\ldots, T-1$ be the trial points and let $\{(B_{t,k},b_{t,k})\}_{k=1}^{|\Xi_t|}$ be the set of realizations used to characterize the random parameters in each stage $t = 2, \ldots, T$. In the same way as defined in Section~\ref{subsec:2stage-partition}, at every stage $t$ we partition $\{(B_{t,k},b_{t,k})\}_{k=1}^{|\Xi_t|}$ into $L_t$ scenario clusters, such that the stage-$t$ \textit{partition} is given by $\N_t = \{\P^\ell_t\}_{\ell=1}^{L_t}$. At stage $t = T$, a \textit{scenario-based subproblem} can be defined for each realization $\xi_{T,k} = (B_{T,k},b_{T,k}) \in \Xi_T$ as follows:
\begin{equation}
\label{eq:eq:scenario-based-t}
\min_{x_T} \left\{c_T^\top x_T \mid A_Tx_T = b_{T,k}-B_{T,k}\bar{x}_{T-1}\right\},
\end{equation}
where a \textit{partition-based subproblem} can be defined similarly to~\eqref{eq:2stage-partition} for each cluster $ \P^\ell_T \in \N_T$ as follows:
\begin{equation}
\label{eq:partition-based-T}
\min_{x_T} \left\{c_T^\top x_T \mid A_Tx_T = \bar{b}_{T,\ell}-\bar{B}_{T,\ell}\bar{x}_{T-1}\right\},
\end{equation}
with $\bar{b}_{T,\ell} := \sum_{k\in \P^\ell_T}b_{T,k}$, and $\bar{B}_{T,\ell} := \sum_{k\in \P^\ell_T}B_{T,k}$.

It can be easily seen that~\eqref{eq:partition-based-T} is constructed by aggregating constraints and variables of~\eqref{eq:eq:scenario-based-t} for all realizations in $\P^\ell_T$, and therefore gives a lower bound for $\sum_{k\in \P^\ell_T}Q_T(\bar{x}_{T-1},\xi_{T,k})$, where $Q_T(\bar{x}_{T-1},\xi_{T,k})$ is the optimal objective value of the scenario-based problem~\eqref{eq:eq:scenario-based-t}. Treating each component $\P^\ell_T$ as a scenario, a coarse cut $\beta^{j\top}_T x_{T-1}+\alpha^j_T \leq r_T$ can be generated in the same way a \textit{standard Benders cut} is generated (see the derivations after equation~\eqref{eq:cpT}), using the corresponding optimal dual solutions of~\eqref{eq:partition-based-T} for each component $\P^\ell_T$ of $\N_T$.

When the coarse cuts do not yield any cut violation at the trial point $\bar{x}_{T-1}$ with respect to the current relaxation for the cost-to-go function, i.e., $\check\Q_T(\bar{x}_{T-1}) \geq \beta^{j\top}_T \bar{x}_{T-1}+\alpha^j_T$, the partition $\N_T$ can be refined by solving the subproblem~\eqref{eq:eq:scenario-based-t} for each $k = 1,2,\ldots, |\Xi_T|$, where the corresponding optimal dual vectors will guide the partition refinements. In this paper we will adopt a very simple refinement strategy by separating two scenarios from a cluster $\P_t$ of the partition $\N_t$ whenever their euclidian distance is sufficiently large. Several other strategies can be considered such as those described in \cite{song2017} for the two-stage case. 

Now consider any other stage $t \in \{2,3,\ldots, T-1\}$, the scenario-based problem now involves the cutting plane approximation for $\Q_{t+1}(x_t)$ which is given by 
\[
\check{\Q}_{t+1}(x_t)=\max_{j\in J_{t+1}} \{\alpha^j_{t+1}+\beta^j_{t+1} x_t\} \  (c.f. \eqref{cpt-1}-\eqref{cutcoeff}). 
\]
Then the scenario-based subproblem for each realization $\xi_{t,k} = (B_{t,k},b_{t,k}) \in \Xi_t$ is given by:
\begin{align}
\label{eq:scenario-based-t}
\underline{Q}_{t,k}(\bar{x}_{t-1}) :=
\left\{
\begin{array}{llll}
 \min \ & c_t^\top x_t+r_{t+1} & \\
\text{s.t. } & A_tx_t= b_{t,k}-B_{t,k}\bar{x}_{t-1}  & (\lambda_k) \\
 & r_{t+1} - \beta^j_{t+1}x_t \geq \alpha^j_t, j\in J_{t+1} \quad & (\pi_{k,j}),
 \end{array}\right.
\end{align}
whereas the partition-based subproblem for each component $\P^\ell_t$ of the partition $\N_t$, for $\ell = 1,2,\ldots, L_t$, can be defined as follows:
\begin{align}
\label{eq:partition-based-t}
\underline{Q}_{t,\ell}(\bar{x}_{t-1}) = 
\left\{ \begin{array}{llll}
\min \ & c_t^\top x_t+r_{t+1} \\
\text{s.t. } & A_t x_t = \bar{b}_{t,\ell}-\bar{B}_{t,\ell}\bar{x}_{t-1}  \\
 & r_{t+1} - \beta^j_{t+1}x_t \geq |\P^\ell_t|\alpha^j_t, j\in J_{t+1},
 \end{array}\right.
\end{align}
which again can be obtained by aggregating variables and constraints of scenario-based subproblems~\eqref{eq:scenario-based-t} for all $k\in \P^\ell_t$.

Similar to the case when $t = T$, the partition refinement is guided by the optimal dual vectors $\{(\lambda_k,\pi_k)\}_{k\in \P^\ell_t}$, i.e., scenarios are put together if the corresponding $(\lambda_k,\pi_k)$ are identical (or "close enough" in euclidean distance). It can be seen from \cite[Theorem 2.5]{song2015} that this partition refinement rule guarantees that, after refining $\N_t$ into $\N'_t = \{\P^\ell_t\}_{\ell=1}^{L'}$,
\[
\sum_{\ell = 1}^{L'_t} \underline{Q}_{t,\ell}(\bar{x}_{t-1}) = \sum_{k=1}^{|\Xi_t|}\underline{Q}_{t,k}(\bar{x}_{t-1}).
\]
Applying the idea of scenario partitions is likely to speed up the process of approximating the cost-to-go functions, particularly in the early stages, where exact information is not important for generating initial cutting-plane relaxations. This has been validated in our experiment results shown in Section~\ref{sec:results}.

\subsection{Adaptive Partition-based SDDP for Multistage Stochastic Linear Programs}
\label{subsec:PartitionSddp}
As thoroughly discussed in \cite{song2017,song2015}, the key reason for the efficiency of the adaptive partition-based strategy in the context of the 2SLP is that, it makes the cut generation effort during the solution procedure adaptive to the solution progress. To demonstrate this further, since most decomposition algorithms used for solving SPs usually rely on sequences of candidate solutions (trial points in the sequel) $\bar x_1=\bar x_1(\xi)$ for generating cutting plane approximations to $\Q_{2}(\cdot)$, it is intuitively clear that generating such approximations to $\Q_{2}(\cdot)$ with high precision for early iterates, likely far from an optimal solution, is surely a wasteful use of computing power --- accuracy will need to be integrated adaptively as the quality of those candidate solutions $\bar x_1=\bar x_1(\xi)$ improves. Ultimately, this procedure not only builds an approximation for $\Q_{2}(\cdot)$ in a computationally efficient manner but also yields a \textit{sufficient} partition $\N$ whose size is (often) smaller than the original number of realizations $|\Xi|$.

Therefore, it is very tempting to conclude that this reduction in the size of the problem will naturally mitigate the computational burden of solving MSLP due to the large number of stages in the planning horizon $T$, and the large number of realizations per stage $|\Xi_t|$. Nonetheless, incorporating the idea of scenario partition for the 2SLP to the MSLP is not necessarily straightforward. To see this, let $\T(\N)$ be a coarse tree induced by the sequence of partitions $\N = (\N_2, \dots, \N_T)$. Suppose we have a procedure by which we sample different sample paths $\xi^s = (\xi^s_1, \xi^s_2, \dots, \xi^s_T)$ from the scenario tree $\T ,\; \forall \; s = 1, 2, \dots, |\Xi_1| \times |\Xi_2| \times \dots \times|\Xi_T|$ and in the backward pass we refine the partition $\N_t , \; \forall \; t=1, 2, \dots T-1$ in order to generate a valid cutting plane approximation for the stage $t$ problem defined over $\xi^s$. Denote a partition that is constructed after taking $p$ sample paths by ${\overrightarrow{\N}}^p$. We say a partition ${\overrightarrow{\N}}^p$ is locally sufficient with respect to $p$ sample paths, if it is sufficient for all $\xi^s \in \{\xi^1, \xi^2, \dots, \xi^p\}$ and we denote such partition by ${{\overrightarrow{\N}}^p}^*$. 

It is not difficult to see that while one partition might be locally sufficient to a given sample path $s$, it is not necessarily sufficient for others. Hence, if one wishes to construct such globally sufficient partition, it is inevitable for a coarse tree constructed after $p$ sample paths $\hat{\T}({\overrightarrow{\N}}^p) \to \T$ as $p \to |\Xi_1| \times \ldots\times |\Xi_T|$. Hence, for a fixed sample path $\xi^s = (\xi^s_1, \xi^s_2, \dots , \xi^s_T)$ one can obtain a sequence of partitions $\N^s = (\N^s_2, \dots , \N^s_T)$ which is only sufficient with respect to $\xi^s$. Additionally, in order to proclaim the sufficiency of $\N^s$ to $\xi^s$ we need to decouple problem~\eqref{eq:mlsp} into a sequence of ($T-1$) consecutive 2SLPs. This decoupling scheme, which we discuss in more details in Subsection~\ref{subsubsec:tree_traversal}, was referred to in \cite{abrahamson1983} and \cite{wittrock1984} for deterministic multistage linear programs as \textit{cautious} and \textit{shuffle}, when the scenario tree is traversed forward and backward, respectively. Consequently, and as we shall see, this not only raises the question of how we can incorporate the adaptive partition-based strategy to the SDDP algorithm, but also raises the question of what the best strategy is to traverse scenario tree.

Moreover, as we intend to employ the adaptive partition-based framework presented in \cite{song2017} which relies on different types of cutting-plane approximations to the cost-to-go functions $\Q_{t+1}(\cdot)$ in their level of inexactness; another aspect to consider in our analysis is, for a given tree traversal strategy, how coarse should the cutting-plane approximation be, such that the computational savings acquired by adaptive partition-based strategies via efficient cut generation effort, offsets the inaccuracy inherited in these inexact cuts during solution process. 

Finally, as a byproduct of this work we attempt to construct an algorithm which exploits the structural nature of the underling problem instance. This can be achieved by integrating this accuracy (inaccuracy) in the quality of cutting planes generated by the adaptive partition-based framework to the SDDP algorithm, relative to how much additional information we can gain from being accurate (as compared to being inaccurate) in different stages of the planning horizon. Next, we characterize in details different ingredients of our proposed algorithms.

\subsubsection{Types of cutting planes}
\label{subsubsec:cut_types}
 In accordance with the definitions introduced in \cite{song2017}, our analysis relies on three types of cutting-plane approximations to the cost-to-go function $\Q_{t+1}(\cdot)$. 
\begin{enumerate}
\item \textit{Fine cuts}. Given an iterate $\bar x_t=\bar x_t(\xi_{[t]})$, a cutting-plane approximation for $\Q_{t+1}(\cdot)$ is generated by solving the scenario-based subproblem~\eqref{eq:scenario-based-t} for each scenario $\xi_{t+1}\in \Xi_{t+1}$.

\item \textit{Coarse cuts}. Given an iterate $\bar x_t=\bar x_t(\xi_{[t]})$, a cutting-plane approximation for $\Q_{t+1}(\cdot)$ is generated by solving the partition-based subproblem~\eqref{eq:partition-based-t} for each cluster  $\P^\ell_{t+1} \in \N_{t+1}$, $\ell = 1,2,\ldots, L_{t+1}$.

\item \textit{Semi-coarse cuts}. We define a \textit{semi-coarse} cut as any hybrid between the \textit{Coarse} and \textit{Fine} cuts. That is, given an iterate $\bar x_t=\bar x_t(\xi_{[t]})$, a cutting-plane approximation for $\Q_{t+1}(\cdot)$ is generated by solving the partition-based subproblem~\eqref{eq:partition-based-t} for $\{\P^{\ell'}_{t+1}\}_{\forall \; \ell' \in L'}$ where $L' \neq \emptyset$ and $L' \subset L$ and solving the scenario-based subproblem~\eqref{eq:scenario-based-t} for each scenario $\xi_{t+1}\in \Xi_{t+1} \backslash \{\P^{\ell'}_{t+1}\}_{\forall \; \ell' \in L'}$.
\end{enumerate}
We refer the reader to \cite{song2017} for a more thorough discussion on these concepts. 
\subsubsection{Tree traversal strategies}
\label{subsubsec:tree_traversal}
We rely on two tree traversal strategies in different variants of the proposed algorithm, some of which were formally introduced in \cite{abrahamson1983} and \cite{wittrock1984} for deterministic multistage linear programs, and further developed in \cite{morton1996} for MSLPs. 

\begin{figure}[htbp]
\begin{center}

\subfigure[\textit{Quick pass}. 
]{
\begin{tikzpicture}
\node [arn_x] (1) at (0,1) {1};
\node [arn_x] (2) at (2.8,1) {2};
\node [arn_x] (3) at (5.6,1) {3};
\node [arn_x] (4) at (8.4,1) {4};
\node [arn_x] (5) at (11.2,1) {5};
\node [arn_w] (1w) at (0,0.5) {};
\node [arn_w] (2w) at (2.8,0.5) {};
\node [arn_w] (3w) at (5.6,0.5) {};
\node [arn_w] (4w) at (8.4,0.5) {};
\node [arn_w] (5w) at (11.2,0.5) {};

\begin{scope}[every node/.style={fill=white, circle, minimum width=0.01cm, minimum height=0.01cm}]
\path [->] (1) edge[fwd] (2);
\path [->] (2) edge[fwd] (3);
\path [->] (3) edge[fwd] (4);
\path [->] (4) edge[fwd] (5);
\path [->] (5w) edge[fwd] (4w);
\path [->] (4w) edge[fwd] (3w);
\path [->] (3w) edge[fwd] (2w);
\path [->] (2w) edge[fwd] (1w);
\end{scope}
\end{tikzpicture}
 \label{fig:TreeTrav-QP}
}

\vspace{1cm}
\subfigure[\textit{Cautious pass}. 
]{
\begin{tikzpicture}

\node [arn_x] (1) at (0,1) {1};
\node [arn_x] (2) at (2.8,1) {2};
\node [arn_x] (3) at (5.6,1) {3};
\node [arn_x] (4) at (8.4,1) {4};
\node [arn_x] (5) at (11.2,1) {5};
\node [arn_w] (1w) at (0,0) {};
\node [arn_w] (2w) at (2.8,0) {};
\node [arn_w] (3w) at (5.6,0) {};
\node [arn_w] (4w) at (8.4,0) {};
\node [arn_w] (5w) at (11.2,0) {};

\begin{scope}[every node/.style={fill=white, circle, minimum width=0.01cm, minimum height=0.01cm}]
\path [->] (1) edge[fwd] (2);
\path [->] (2) edge[fwd] (3);
\path [->] (3) edge[fwd] (4);
\path [->] (4) edge[fwd] (5);
\path [->] (1w) edge[fwd,  bend left=50] (2w);
\path [->] (2w) edge[fwd,  bend left=50] (3w);
\path [->] (3w) edge[fwd,  bend left=50] (4w);
\path [->] (4w) edge[fwd,  bend left=50] (5w);
\path [->] (5w) edge[fwd,  bend left=50] (4w);
\path [->] (4w) edge[fwd,  bend left=50] (3w);
\path [->] (3w) edge[fwd,  bend left=50] (2w);
\path [->] (2w) edge[fwd,  bend left=50] (1w);

\end{scope}
\end{tikzpicture}
 \label{fig:TreeTrav-CP}
}

\subfigure{
\begin{tikzpicture}
\matrix [draw,below left] at (current bounding box.south) {
\node [arn_vsw] (0bd) at (-2,0) {};
\node [arn_vsw,label=right:Backward pass] (00bd) at (0.25,0) {};
\path [->] (00bd) edge[fwd, line width = 0.5 mm] (0bd);\\
\node [arn_vsw] (1bd) at (-2,0) {};
\node [arn_vsw,label=right: Forward pass] (11bd) at (0.25,0) {};
\path [->] (1bd) edge[fwd, line width = 0.5 mm] (11bd);\\
\node [arn_vsw] (6bd) at (-2,0) {};
\node [arn_vsw,label=right:Inner forward step] (66bd) at (0.25,0) {};
\path [->] (6bd) edge[fwd, bend left=40,  line width = 0.5 mm] (66bd);\\
\node [arn_vsw] (7bd) at (-2,0) {};
\node [arn_vsw,label=right:Inner backward step] (77bd) at (0.25,0) {};
\path [->] (77bd) edge[fwd, bend left=40, line width = 0.5 mm] (7bd);\\
};
\end{tikzpicture}
}
\end{center}
\caption{Tree traversal strategies.}
\label{fig:TreeTrav}
\end{figure}
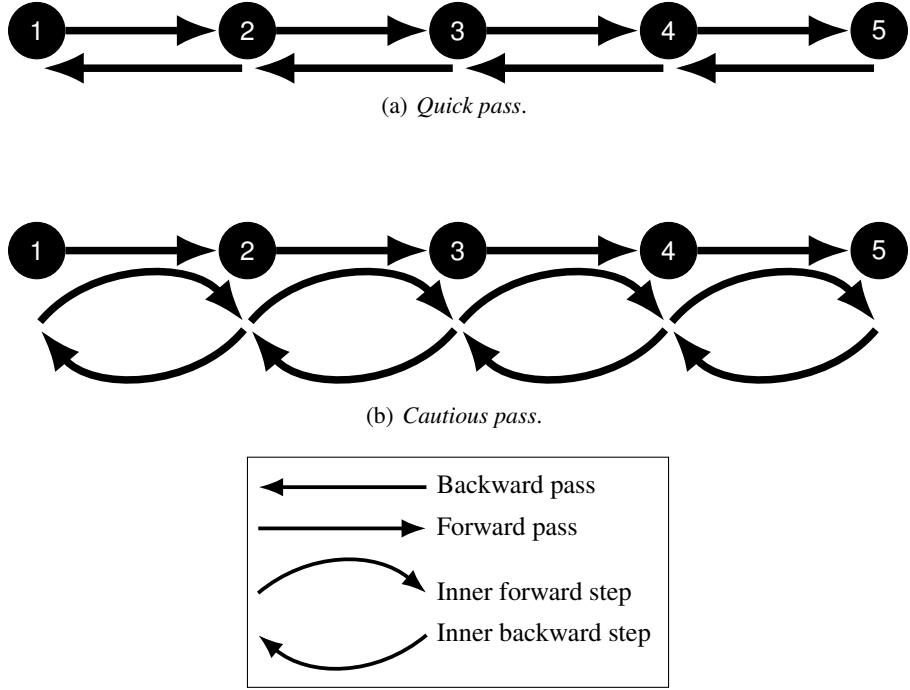

\paragraph{Quick pass (QP).} Under this scheme, the policy (which is induced by the current approximation for the cost-to-go functions $\Q_{t+1}(\cdot)$) evaluation process iteratively passes candidate solution $\bar x_t=\bar x_t(\xi_{[t]})$ down the scenario tree stage by stage ($t= 1 \to t=2 \to \dots t=T$) and cuts (if any exists) are passed directly back up the tree  ($t= T \to t=T-1 \to \dots t=1$) with no intermediate change of direction between any two consecutive stages $t$ and $t-1$. We care to mention that \textit{quick pass} is in fact, the most commonly used strategy in the SDDP algorithm. QP is illustrated in Figure~\ref{fig:TreeTrav-QP}.
\paragraph{Cautious pass (CP).} Under this scheme, the iterative policy evaluation process never goes back up the tree unless all cuts that would be passed back from stage $t+1$ to stage $t, \; \forall \; t=T-1, T-2, \dots 1$, are redundant, i.e., the evaluation process will maintain intermediate changes of directions between consecutive stages continuously. Note that, an intermediate change of direction going forward entails, updating the candidate solution $\bar x_t=\bar x_t(\xi_{[t]})$ (if possible) based on the updated cutting-plane approximation $\check \Q'_{t+1}(\cdot)$. We refer to such intermediate changes of direction as, \textit{inner forward step} and \textit{inner backward step} when the change of direction is forward and backward, respectively. CP is illustrated in Figure~\ref{fig:TreeTrav-CP}.

\subsubsection{Refinements strategies}
\label{subsubsec:refinements}
The specialization of the proposed adaptive partition-based SDDP algorithm depends on how a partition refinement rule is integrated with the SDDP algorithm. The particular partition refinement rule that we used in this paper is motivated by Song and Luedtke~\cite{song2015} under the name of the absolute rule, where partition refinement is guided by the optimal dual vectors $\{(\lambda_{\xi_t},\pi_{\xi_t})\}_{\xi_t \in \P^\ell_t}$. That is, any two realizations $\xi^i_t$ and $\xi^j_t \; \forall \; i \neq j$ from a cluster $\P^\ell_t$ of the partition $\N_t$ will be split up whenever their euclidean distance is sufficiently large, i.e., $\norm{(\lambda_{\xi^i_t},\pi_{\xi^i_t})-(\lambda_{\xi^j_t},\pi_{\xi^j_t})} > \epsilon$, where $\epsilon>0$ is a user-specified tolerance parameter. We propose the following two schemes for integrating the partition refinements to the SDDP algorithm:
\begin{enumerate}
\item \textit{Refinement within SDDP}. Given the original scenario tree $\T$ and a coarse tree $\hat{\T}(\N)$ induced by a sequence of partitions $\N = (\N_2, \dots , \N_T)$, sample from the original scenario tree $\T$ going forward and continuously attempt to refine the partition $\N_t$, $\forall \; t=2, \dots \T$ while generating cuts during the backward pass of the SDDP algorithm. 

\item \textit{Refinement outside SDDP}. Perform the SDDP algorithm on a coarse tree $\hat{\T}(\N)$, i.e., both sampling and cut generation are from the coarse tree only, without any attempt to refine any partition $\N_t$, $\forall \; t=2, \dots T$ during the backward pass, and refine the sequence of partitions $\N$ in a separate refinement step whenever it is needed. This separate step is introduced in which the refinement is performed locally on sample path(s) $\xi^s = (\xi_1^s, \dots , \xi_{T-1}^s)$ that is sampled from the original scenario tree $\T$. 
\end{enumerate} 
The main feature of the \textit{Refinement within SDDP} strategy is its flexibility to add cuts from both the coarse tree $\hat{\T}(\N)$ and original tree $\T$ at each step of the SDDP algorithm during the solution process. However, one concern with this approach is that, at any stage $t$ whenever the coarse cut $\bar q_t(x_{t-1})$ (generated by aggregating the variables and constraints according to $\N_t$) does not yield any violation to the current iterate, the process will always attempt to refine the partition $\N_t, \; \forall \; t=T, \dots 2$. Hence, the partition size can grow very quickly as we step into different sample paths -- which makes it difficult to effectively exploit the reduction in the size of the problem.

On the other hand, while the \textit{Refinement outside SDDP} strategy is more restricted in the sense that, it generates cuts from the coarse tree only during the execution of the SDDP algorithm, it capitalizes on the merits of the adaptive partition-based framework, since the SDDP algorithm is now performed on a smaller scenario tree. Nonetheless, the main concern with the \textit{Refinement outside SDDP} strategy is the difficulty of identifying a criterion by which we can claim the insufficiency of the current sequence of partitions $\N$, which would suggest its refinement.

\section{Implementation Details}
\label{sec:ImplmntDetails}
In this section we describe in more detail the proposed integrated adaptive partition-based SDDP algorithms. As previously noted, it is almost impossible to construct a sufficient partition which accommodates all possible sample paths in the scenario tree. However, depending on the refinement rule, one can easily have some control over the size of the coarse tree $\hat{\T}({\overrightarrow{\N}}^p)$ obtained after $p$ sample paths. This refinement rule can be driven by the manner in which the scenario tree is traversed based on different strategies introduced in Subsection~\ref{subsubsec:tree_traversal}. On one hand, employing a CP is clearly a more rigorous strategy in refining the sequence of partitions ${\overrightarrow{\N}}^p$, since it attempts to create a sufficient sequence of partitions ${{\overrightarrow{\N}}^p}^*$ by solving a sequence of $T-1$ 2SLP problems (defined by \nth{1}-stage $=t$ and \nth{2}-stage $=t+1, \; \forall \; t=1, 2, \dots T-1$) to optimality. On the other hand, a QP is a more of a lenient strategy to refine ${\overrightarrow{\N}}^p$ since it is only attempting to refine the partition ${\overrightarrow{\N}}^p$ without any intention to achieve local sufficiency.

While it is instructive to construct a locally sufficient sequence partitions ${{\overrightarrow{\N}}^p}^*$ from an optimality point of view; as we shall see in Section~\ref{sec:results}, solving a sequence of $T-1$ 2SLP problems up to optimality could be expensive, and perhaps not worth doing, especially at early iterates when the candidate solution $\bar x_t=\bar x_t(\xi_{[t]})$ is likely to be far from optimal. To that end, as mentioned in Subsection~\ref{subsubsec:refinements} we consider two refinement schemes for integrating the adaptive partition-based approach to the SDDP algorithm, namely Refinement within SDDP, and Refinement outside SDDP, which will be described in details in the following two subsections, respectively.

\subsection{Refinement outside SDDP}
\label{subsec:RefineOut}
We present two strategies of performing partition refinement outside SDDP:
\begin{itemize}
\item \textit{Adaptive partition-enabled preprocessing for SDDP (APEP-SDDP)}: the adaptive partition-based approach is used only as a preprocessing step for the SDDP algorithm, which is applied on the original scenario tree. A coarse scenario tree is iteratively refined during the preprocessing step. Once the size of the coarse scenario tree is sufficiently large, the preprocessing step is terminated and the SDDP algorithm will be used for the original scenario tree for the rest of the procedure.
\item \textit{SDDP on iteratively refined coarse scenario trees (ITER-SDDP)}: apply the SDDP algorithm for each coarse scenario tree, which is iteratively refined throughout the solution process. 
\end{itemize}
 
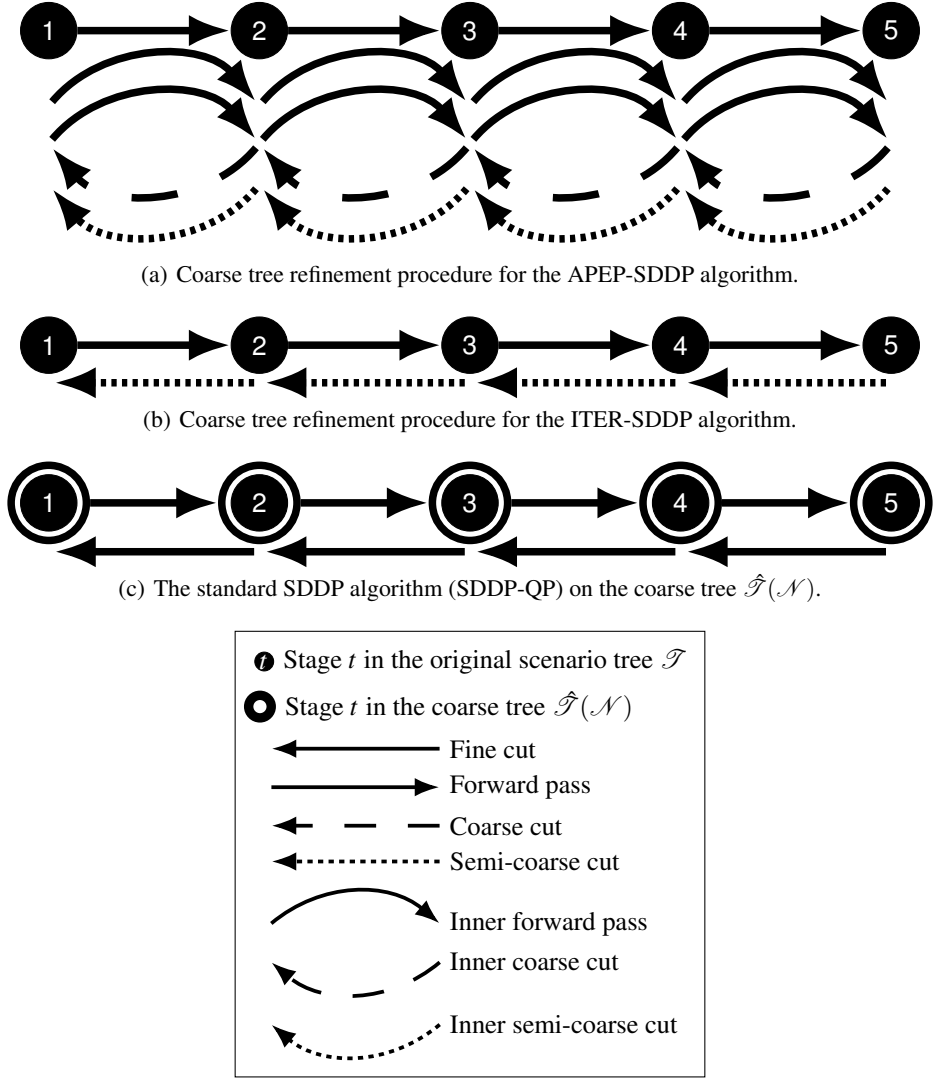
\begin{figure}[htbp]
\begin{center}
\subfigure[Coarse tree refinement procedure for the APEP-SDDP algorithm.]{
\begin{tikzpicture}
\node [arn_x] (1) at (0,1) {1};
\node [arn_x] (2) at (2.8,1) {2};
\node [arn_x] (3) at (5.6,1) {3};
\node [arn_x] (4) at (8.4,1) {4};
\node [arn_x] (5) at (11.2,1) {5};
\node [arn_w] (1w) at (0,0) {};
\node [arn_w] (2w) at (2.8,0) {};
\node [arn_w] (3w) at (5.6,0) {};
\node [arn_w] (4w) at (8.4,0) {};
\node [arn_w] (5w) at (11.2,0) {};
\node [arn_w] (1wc) at (0,-0.5) {};
\node [arn_w] (2wc) at (2.8,-0.5) {};
\node [arn_w] (3wc) at (5.6,-0.5) {};
\node [arn_w] (4wc) at (8.4,-0.5) {};
\node [arn_w] (5wc) at (11.2,-0.5) {};

\node [arn_w] (1w1) at (0,-1.05) {};
\node [arn_w] (2w1) at (2.8,-1.05) {};
\node [arn_w] (3w1) at (5.6,-1.05) {};
\node [arn_w] (4w1) at (8.4,-1.05) {};
\node [arn_w] (5w1) at (11.2,-1.05) {};
\node [arn_w] (1wc1) at (0,-0.505) {};
\node [arn_w] (2wc1) at (2.8,-0.505) {};
\node [arn_w] (3wc1) at (5.6,-0.505) {};
\node [arn_w] (4wc1) at (8.4,-0.505) {};
\node [arn_w] (5wc1) at (11.2,-0.505) {};

\begin{scope}[every node/.style={fill=white, circle, minimum width=0.01cm, minimum height=0.01cm}]
\path [->] (1) edge[fwd] (2);
\path [->] (2) edge[fwd] (3);
\path [->] (3) edge[fwd] (4);
\path [->] (4) edge[fwd] (5);
\path [->] (5w1) edge[Fwdds, bend left=50] (4w1);
\path [->] (4w1) edge[Fwdds, bend left=50] (3w1);
\path [->] (3w1) edge[Fwdds, bend left=50] (2w1);
\path [->] (2w1) edge[Fwdds, bend left=50] (1w1);
\path [->] (5wc1) edge[fwdds, bend left=50] (4wc1);
\path [->] (4wc1) edge[fwdds, bend left=50] (3wc1);
\path [->] (3wc1) edge[fwdds, bend left=50] (2wc1);
\path [->] (2wc1) edge[fwdds, bend left=50] (1wc1);

\path [->] (1w) edge[fwd, bend left=50] (2w);
\path [->] (2w) edge[fwd, bend left=50] (3w);
\path [->] (3w) edge[fwd, bend left=50] (4w);
\path [->] (4w) edge[fwd, bend left=50] (5w);
\path [->] (1wc) edge[fwd, bend left=50] (2wc);
\path [->] (2wc) edge[fwd, bend left=50] (3wc);
\path [->] (3wc) edge[fwd, bend left=50] (4wc);
\path [->] (4wc) edge[fwd, bend left=50] (5wc);

\end{scope}
\end{tikzpicture}
\label{fig:APEP-SDDP}
}

\subfigure[Coarse tree refinement procedure for the ITER-SDDP algorithm.]{
\begin{tikzpicture}
\node [arn_x] (1) at (0,1) {1};
\node [arn_x] (2) at (2.8,1) {2};
\node [arn_x] (3) at (5.6,1) {3};
\node [arn_x] (4) at (8.4,1) {4};
\node [arn_x] (5) at (11.2,1) {5};
\node [arn_w] (1w) at (0,0.5) {};
\node [arn_w] (2w) at (2.8,0.5) {};
\node [arn_w] (3w) at (5.6,0.5) {};
\node [arn_w] (4w) at (8.4,0.5) {};
\node [arn_w] (5w) at (11.2,0.5) {};

\begin{scope}[every node/.style={fill=white, circle, minimum width=0.01cm, minimum height=0.01cm}]
\path [->] (1) edge[fwd] (2);
\path [->] (2) edge[fwd] (3);
\path [->] (3) edge[fwd] (4);
\path [->] (4) edge[fwd] (5);
\path [->] (5w) edge[Fwdds] (4w);
\path [->] (4w) edge[Fwdds] (3w);
\path [->] (3w) edge[Fwdds] (2w);
\path [->] (2w) edge[Fwdds] (1w);
\end{scope}
\end{tikzpicture}
\label{fig:ITER-SDDP}
}

\subfigure[The standard SDDP algorithm (SDDP-QP) on the coarse tree $\hat{\T}(\N)$.]{
\begin{tikzpicture}
\node [sqr_x] (11) at (0,1) {1};
\node [sqr_x] (22) at (2.8,1) {2};
\node [sqr_x] (33) at (5.6,1) {3};
\node [sqr_x] (44) at (8.4,1) {4};
\node [sqr_x] (55) at (11.2,1) {5};
\node [arn_x] (1) at (0,1) {1};
\node [arn_x] (2) at (2.8,1) {2};
\node [arn_x] (3) at (5.6,1) {3};
\node [arn_x] (4) at (8.4,1) {4};
\node [arn_x] (5) at (11.2,1) {5};
\node [arn_w] (1w) at (0,0.35) {};
\node [arn_w] (2w) at (2.8,0.35) {};
\node [arn_w] (3w) at (5.6,0.35) {};
\node [arn_w] (4w) at (8.4,0.35) {};
\node [arn_w] (5w) at (11.2,0.35) {};

\begin{scope}[every node/.style={fill=white, circle, minimum width=0.01cm, minimum height=0.01cm}]
\path [->] (11) edge[fwd] (22);
\path [->] (22) edge[fwd] (33);
\path [->] (33) edge[fwd] (44);
\path [->] (44) edge[fwd] (55);
\path [->] (5w) edge[fwd] (4w);
\path [->] (4w) edge[fwd] (3w);
\path [->] (3w) edge[fwd] (2w);
\path [->] (2w) edge[fwd] (1w);
\end{scope}
\end{tikzpicture}
\label{fig:PartsSDDPalgrthms-coarseSDDP}
}

\subfigure{
\begin{tikzpicture}
\matrix [draw,below left] at (current bounding box.south) {
\node [arn_xx,label=right: Stage $t$ in the original scenario tree $\T$]at (-2,0) {$t$}; \\
\node [sqr_xx,label=right:Stage $t$ in the coarse tree $\hat{\T}(\N)$ ] at (-2,0) {}; \\
\node [arn_vsw] (0bd) at (-2,0) {};
\node [arn_vsw,label=right:Fine cut] (00bd) at (0.25,0) {};
\path [->] (00bd) edge[fwd, line width = 0.5 mm] (0bd);\\
\node [arn_vsw] (1bd) at (-2,0) {};
\node [arn_vsw,label=right: Forward pass] (11bd) at (0.25,0) {};
\path [->] (1bd) edge[fwd, line width = 0.5 mm] (11bd);\\
\node [arn_vsw] (3bd) at (-2,0) {};
\node [arn_vsw,label=right:Coarse cut] (33bd) at (0.25,0) {};
\path [->] (33bd) edge[fwdds,  line width = 0.5 mm, dash pattern = on 12pt] (3bd);\\
\node [arn_vsw] (4bd) at (-2,0) {};
\node [arn_vsw,label=right: Semi-coarse cut] (44bd) at (0.25,0) {};
\path [->] (44bd) edge[Fwdds,  line width = 0.5 mm] (4bd);\\
\node [arn_vsw] (6bd) at (-2,0) {};
\node [arn_vsw,label=right:Inner forward pass] (66bd) at (0.25,0) {};
\path [->] (6bd) edge[fwd, bend left=40,  line width = 0.5 mm] (66bd);\\
\node [arn_vsw] (7bd) at (-2,0) {};
\node [arn_vsw,label=right:Inner coarse cut] (77bd) at (0.25,0) {};
\path [->] (77bd) edge[fwdds, bend left=40, line width = 0.5 mm,  dash pattern = on 12pt] (7bd);\\
\node [arn_vsw,label=right:Inner semi-coarse cut] (77bd) at (0.25,0) {};
\path [->] (77bd) edge[Fwdds, bend left=40, line width = 0.5 mm] (7bd);\\
};
\end{tikzpicture}
}
\end{center}
\caption{Illustration of tree traversal strategies and various types of cutting planes used in the \textit{Refinement outside SDDP} algorithms.}
\label{fig:RO-SDDP}
\end{figure}

In both strategies, we determine whether the sequence of partitions $\N$ needs to be refined after applying $j$ iterations of the standard SDDP algorithm on the coarse tree $\hat{\T}({\N})$, by keeping track of the lower bound progress for the last $n$ consecutive iterations.  If $\underline{z}^{j} - \underline{z}^{j-n} \leq \epsilon$, where $\epsilon > 0$ is a given tolerance, it indicates that the lower bound does not improve by more than $\epsilon$ over the last $n$ iterations. Both APEP-SDDP and ITER-SDDP will then refine the sequence of partitions $\N$ to $\N'$, and perform the SDDP algorithm on the refined tree $\hat{\T}({\N'})$, except that if the coarse scenario tree after the refinement is large enough, the APEP-SDDP algorithm will revert back to the original scenario tree and perform the standard SDDP algorithm on it afterwards. Figure~\ref{fig:RO-SDDP} illustrates these two algorithms for \textit{Refinement outside SDDP}. 

\subsubsection{The APEP-SDDP algorithm}
\label{subsubsec:APEP-SDDP}
The primary concern with the APEP-SDDP algorithm is the difficulty to identify an adequate criterion for terminating the preprocessing step without jeopardizing the efficacy of the algorithm. Our implementations adopt a heuristic criterion for terminating the preprocessing step based on the size of the coarse tree $\hat{\T}({\N})$. More specifically, we define the size of the coarse tree as $|\hat{\T}({\N})| =  \frac{\sum_{t=2}^{T} {\N}_t/|\Xi_t|}{T-1}$ and let $\nu \in (0,1)$ be a user-specified parameter for terminiating the preprocessing step. When $|\hat{\T}({\overrightarrow{\N}}^p)| > \nu$, we terminate the preprocessing and revert to solving the original scenario tree $\T$ using the standard SDDP algorithm afterwards. Alternatively, one might choose to define the termination criterion as a fixed fraction of the computational budget (time-limit) and set $\nu$ accordingly.  Finally, we use CP as the tree traversal strategy when refining $\N$. The motivation of this choice is that, using CP to traverse $\T$ is more likely to result in a larger $|\hat{\T}({\N})|$, making the preprocessing step to be terminated more quickly. We have found that this choice yielded improved performance in our numerical experiments. The APEP-SDDP algorithm is summarized in Algorithm~\ref{alg:APEP-SDDP}. 
\begin{algorithm}[htbp]
\caption{The APEP-SDDP algorithm.}
\label{alg:APEP-SDDP}
\begin{algorithmic}[0]
\State \textbf{STEP 0:} Initialization. 
\begin{enumerate}
\item  Let $p=0,\; \underline{z}^{p}:= -\infty,\; \check \Q^{p}_t(\cdot) := -\infty, \; \forall \; t=2, \dots, T$. 
\item Define an initial sequence of partitions ${\overrightarrow{\N}}^p=({\overrightarrow{\N}}^p_2,\ldots,{\overrightarrow{\N}}^p_T)$ and the corresponding coarse tree  $\hat{\T}({\overrightarrow{\N}}^p)$.
\item Define a parameter $n > p$.
\item Choose a tolerance $\epsilon > 0$ and a preprocessing termination parameter $\nu \in (0,1)$.
\end{enumerate}
\State \textbf{STEP 1:} 
\begin{enumerate}
\item Increment $p \gets  p+1$.
\item If $|\hat{\T}({\overrightarrow{\N}}^p)| > \nu$, go to STEP 3. Otherwise, define an MSLP on $\hat{\T}({\overrightarrow{\N}}^p$) and the cutting plane approximations $\check \Q^{p}_t(\cdot) \; \forall \; t=2, \dots, T$ as follows:
\begin{enumerate}[label=(\roman*)]
\item Initialize $j_p > n$ and set $\underline{z}^{j_p} = \underline{z}^{p}$ and $\check \Q^{j_p}_t(\cdot) = \check \Q^{p}_t(\cdot),\; \forall \; t=2, \dots T$.
\item \textbf{while} {\textit{true}} \textbf{do}
{
\begin{enumerate}[label=\alph*.]
\item  Increment $j_p \gets j_p+1$. 
\item Apply the forward and backward step of the standard SDDP algorithm to improve the cutting plane approximations and update $\check \Q^{j_p}_t(\cdot),\; \forall \; t=2, \dots, T$ and $\underline{z}^{j_p}$.
\item \textbf{if} $\underline{z}^{j_p} - \underline{z}^{j_p-n} \leq \epsilon$
\begin{itemize}
\item[] set  $\underline{z}^{p}= \underline{z}^{j_p}$, and $\check \Q^{p}_t(\cdot) = \check \Q^{j_p}_t(\cdot),\; \forall \; t=2, \dots T$; BREAK.
\end{itemize}
\textbf{end} 
\end{enumerate}}
\item[] \textbf{end}
\end{enumerate}
\end{enumerate}
\State \textbf{STEP 2:}
\begin{enumerate}
\item Choose a sample path $\xi^p = (\xi^p_1, \dots , \xi^p_{T-1})$. 
\item Refine ${\overrightarrow{\N}}^p$ by a \textbf{CP} and the adaptive partition-based approach \cite{song2017,song2015} over the sample path $\xi^p$ to construct $\hat{\T}({\overrightarrow{\N}}^{p+1})$.
\item Update $\underline{z}^{p}$ and $ \check \Q^{p}_t(\cdot)$.
\item Go to STEP 1.  
\end{enumerate}
\State \textbf{STEP 3:}
\begin{enumerate}
\item Define an MSLP on $\T$ and initialize the approximate cost-to-go functions using $ \check \Q^{p}_t(\cdot)$.
\item Solve it using the standard SDDP algorithm and stop upon a termination criterion (such as a time limit).
\end{enumerate}
\end{algorithmic}
\end{algorithm}

\subsubsection{The ITER-SDDP algorithm}
\label{subsubsec:ITER-SDDP}
As previously noted, the ITER-SDDP algorithm extends the APEP-SDDP algorithm in the sense that, here, there is no threshold by which we stop generating cutting planes from the coarse tree $\hat{\T}({\overrightarrow{\N}}^p)$ and revert back to using the standard SDDP algorithm on the original scenario tree $\T$. In other words, the ITER-SDDP algorithm will continuously attempt to refine $\hat{\T}({\overrightarrow{\N}}^p)$, for every iteration $p$ whenever the lower bound $\underline{z}$ does not make significant progress. In this case, the cut generation from the coarse tree is no longer a "preprocessing" step, but integrated into the entire solution procedure. From the algorithmic perspective, this can be easily achieved by setting $\nu = 1$ in Algorithm~\ref{alg:APEP-SDDP}, in which case the algorithm will never reach \textbf{STEP 3}. 

In addition, unlike the CP traversal strategy employed in \textbf{STEP 2} (the refinement step) in Algorithm~\ref{alg:APEP-SDDP}, in the ITER-SDDP algorithm we use QP as the tree traversal strategy. As previously noted, using QP which is a more lenient strategy than CP, will allow $|\hat{\T}({\N})|$ to be kept small during the solution process. 

\subsection{Refinement within SDDP}
\label{subsec:RefineWithin}
Unlike the \textit{Refinement outside SDDP} strategy introduced in Section~\ref{subsec:RefineOut}, where the partition refinement is performed separately from the SDDP algorithm, whenever the lower bound progress is not significant during the SDDP algorithm. The \textit{Refinement within SDDP} strategy is one in which the partition refinement is attempted at every \textit{backward pass} of the SDDP algorithm if necessary, irrespective of the progress in the lower bound. Here, both the partition refinement and the SDDP algorithm are performed based on sample paths generated from the original scenario tree $\T$. In our implementations we consider two variants of the Refinement within SDDP strategy, which will be described in more details next.

\subsubsection{The adaptive partition-based SDDP algorithm with quick cut generation (APQP-SDDP)}
\label{subsubsec:APQP-SDDP}
This algorithm can be seen as the most natural extension of integrating the adaptive partition-based strategies to the backward pass of the standard SDDP algorithm. In the backward pass, at every stage $t = T, \dots, 2$, the APQP-SDDP algorithm starts by attempting to generate a \textit{coarse cut}, and if it succeeds in doing so, the process immediately moves on to stage $t-1$; otherwise, it attempts to generate a \textit{semi-coarse cut} by sequentially going through different clusters $\P^{\ell}_{t} \in {\overrightarrow{\N}}^p_t$, and as soon as this process succeeds in generating a violated cut, it updates $\check \Q^{p}_t(\cdot)$, refines ${\overrightarrow{\N}}^p_t$ and moves on to stage $t-1$. In other words, any violated cut, regardless of its quality, is a sufficient criterion for moving back to stage $t-1$. We summarize the APQP-SDDP algorithm in Algorithm~\ref{alg:APQP-SDDP} and illustrate it in Figure~\ref{fig:PartsSDDPalgrthms-APQP}.

\begin{algorithm}[htbp]
\caption{The adaptive partition-based SDDP algorithm with quick cut generation (APQP-SDDP).}
\label{alg:APQP-SDDP}
\begin{algorithmic}[0]
\State \textbf{STEP 0:} Initialization.
\begin{enumerate}
\item Let $p=0, \; \check \Q^{p}_t(\cdot) := -\infty, \; \epsilon \geq 0.$
\item Define an initial sequence of partitions ${\overrightarrow{\N}}^p=({\overrightarrow{\N}}^p_1,\ldots,{\overrightarrow{\N}}^p_T)$, the corresponding coarse tree $\hat{\T}({\overrightarrow{\N}}^p)$.
\item Initialize the corresponding necessary parameters of the SDDP algorithm. 
\end{enumerate} 
\State \textbf{STEP 1:} Increment $p \gets p+1$ and choose a sample path $\xi^p = (\xi^p_1, \dots , \xi^p_T)$.
\State \textbf{STEP 2:} Implement the SDDP \textit{forward step} over $\xi^p$ to obtain $\bar x_t=\bar x_t(\xi_{[t]}^p),\; \forall t=1, \dots T-1$.
\State \textbf{STEP 3:} Implement the adaptive partition-based SDDP \textit{Backward-step} as follows:
\begin{enumerate}
\item For $t= T-1, \dots, 1$
\begin{enumerate}
\item Solve the partition-based subproblem~\eqref{eq:partition-based-t} for each cluster  $\P^\ell_{t+1} \in {\overrightarrow{\N}}^p_{t+1}$, $\ell = 1,2,\ldots, L_{t+1}$ to improve the cutting plane approximation for $\check\Q_{t+1}(\bar x_t)$.
\item If $\bar\Q_{t+1}(\bar x_t) - \check \Q^{p}_{t+1}(\bar x_t) > \epsilon$, update $\check \Q_{t+1}(x_t)$ and go to $t-1$. Else, go to (c).
\item Compute a semi-coarse cut as described in~\eqref{subsubsec:cut_types} and update $\check \Q_{t+1}(x_t)$.
\end{enumerate}
\end{enumerate}
\vspace{0.1cm}
\State \textbf{STEP 4:} If termination criterion is achieved, STOP. Else, go to STEP 1. 
\end{algorithmic}
\end{algorithm}

\subsubsection{The adaptive partition-based SDDP algorithm with cautious cut generation (APCP-SDDP)}
\label{subsubsec:APCP-SDDP}
Instead of using the existence of any cut violation regardless of its quality as the criterion for moving back to the previous stage, as what is used in APQP-SDDP, the APCP-SDDP algorithm attempts to generate all possible cuts of various types by employing inner forward and backward steps just as SDDP-CP. In our implementation of APCP-SDDP, we follow the same procedure described in Algorithm~\ref{alg:APQP-SDDP} except that here, we modify \textbf{STEP 3} to:
\begin{itemize}
\item[$\bullet$] $\forall \; t=T-1, T-2, \dots 1$, implement the adaptive partition-based approach of \cite{song2017} to solve a 2SLP problem over $\xi^p_t$ which is defined by $t$ as \nth{1}-stage and $t+1$ as the \nth{2}-stage.
\end{itemize}
Figure~\ref{fig:PartsSDDPalgrthms-APCP} illustrates the APCP-SDDP algorithm.

\subsection{Adaptive partition-based SDDP with structured cut generation policies}
\label{subsubsec:SPAP-SDDP}
In the previous two subsections,~\ref{subsec:RefineOut} and~\ref{subsec:RefineWithin}, we provide algorithms in which the same type of cut generation strategies are implemented uniformly across all stages $t=2, \dots, T$. In this subsection, we investigate a framework which incorporates different types of cut generation strategies in different stages. To that end, we classify different stages in the planning horizon into two different categories:
\begin{enumerate}
\item Stages of "less-importance", where we use \textit{coarse} and \textit{semi-coarse} cuts.
\item Stages of "more-importance", where we use \textit{fine} cuts only.
\end{enumerate}
The motivation for associating \textit{coarse} and \textit{semi-coarse} cuts with stages considered to be less important and \textit{fine} cuts with stages considered to be more important is as follows. In the MSLP setting, even if we simply attempt to generate $\textit{fine cuts}$ without any aggregation, the cutting hyperplane $q_t(x_{t-1}) := \beta^\top_t x_{t-1}+\alpha_t$ generated in stages $t = T-1, \dots, 2$ still might not be a supporting hyperplane to $\Q_t(\cdot)$, since the approximation error propagates as $t = T-1 \to t=2$. Hence, using aggregated information with respect to a partition $\N^{p}_{t}$ might hinder the performance of the algorithm by adding another layer of inaccuracy. In certain situations, we may not be able to afford this additional layer of inaccuracy in stages where we need to be more accurate, in which case we have to use \textit{fine cuts}. On the other hand, in stages where we believe there is little extra information to be gained from \textit{fine} cuts compared to \textit{coarse} cuts, we can save some computational effort by utilizing cutting-plane approximations of any quality (\textit{coarse} and \textit{semi-coarse} cuts). We refer to this algorithm as the SPAP-SDDP algorithm, summarize it in Algorithm~\ref{alg:SPAP-SDDP} and illustrate it in Figure~\ref{fig:PartsSDDPalgrthms-SPAP}.

\begin{algorithm}[htbp]
\caption{The adaptive partition-based SDDP algorithm with structured cut generation policies (SPAP-SDDP).}
\label{alg:SPAP-SDDP}
\begin{algorithmic}[0]
\State \textbf{STEP 0:} Initialization.
\begin{enumerate}
\item Let $p=0, \; \check \Q^{p}_t(\cdot) := -\infty, \; \epsilon \geq 0, \; \Z  \in \mathbb{R}$
\item Classify every stage $t=2, \dots, T$ such that $t$ is "more-important" if $t \in MI := \{t \in \mathbb{R}_{+} |  \bar{z}_t \leq \Z \}$ and $t$ is "less-important" if $t \in LI := \{t \in \mathbb{R}_{+} |  \bar{z}_t > \Z \}$.
\item Define an initial sequence of partitions ${\overrightarrow{\N}}^p=({\overrightarrow{\N}}^p_1,\ldots, {\overrightarrow{\N}}^p_T)$, the corresponding coarse tree $\hat{\T}({\overrightarrow{\N}}^p), \; \forall \; t \in LI $
\item Initialize the corresponding necessary parameters of the SDDP algorithm. 
\end{enumerate} 
\State \textbf{STEP 1:} Increment $p \gets p+1$ and choose a sample path $\xi^p = (\xi^p_1, \dots , \xi^p_T)$.
\State \textbf{STEP 2:} Implement the SDDP \textit{forward step} over $\xi^p$ to obtain $\bar x_t=\bar x_t(\xi_{[t]}),\; \forall t=1, \dots T-1$.
\State \textbf{STEP 3:} Implement the \textit{backward step} as follows:
\begin{enumerate}
\item For $t= T-1, \dots, 1$
\begin{enumerate}
\item If $t \in MI$ compute a \textit{fine cut} by solving the scenario-based subproblem~\eqref{eq:scenario-based-t} for each scenario $\xi_t \in \Xi_t$. Else, go to (b).
\item solve the partition-based subproblem~\eqref{eq:partition-based-t} for each cluster  $\P^\ell_{t+1} \in {\overrightarrow{\N}}^p_{t+1}$, $\ell = 1,2,\ldots, L_{t+1}$ to compute a coarse cut for $\bar\Q_{t+1}(\bar x_t)$.
\item If $\bar\Q_{t+1}(\bar x_t) - \check \Q^{p}_{t+1}(\bar x_t) > \epsilon$, update $\check \Q_{t+1}(x_t)$ and go to $t-1$. Else, go to (d).
\item compute a semi-coarse cut as described in~\eqref{subsubsec:cut_types} and $\check \Q_{t+1}(x_t)$.
\end{enumerate}
\end{enumerate}
\vspace{0.1cm}
\State \textbf{STEP 4:} If termination criterion is achieved, STOP. Else, go to STEP 1. 
\end{algorithmic}
\end{algorithm}

\begin{figure}[htbp]
\begin{center}

\subfigure[The adaptive partition-based SDDP algorithm with quick cut generation (APQP-SDDP).]{
\begin{tikzpicture}
\node [arn_x] (1) at (0,1) {1};
\node [arn_x] (2) at (2.8,1) {2};
\node [arn_x] (3) at (5.6,1) {3};
\node [arn_x] (4) at (8.4,1) {4};
\node [arn_x] (5) at (11.2,1) {5};
\node [arn_w] (1w) at (0,0.5) {};
\node [arn_w] (2w) at (2.8,0.5) {};
\node [arn_w] (3w) at (5.6,0.5) {};
\node [arn_w] (4w) at (8.4,0.5) {};
\node [arn_w] (5w) at (11.2,0.5) {};
\node [arn_w] (1wc) at (0,0) {};
\node [arn_w] (2wc) at (2.8,0) {};
\node [arn_w] (3wc) at (5.6,0) {};
\node [arn_w] (4wc) at (8.4,0) {};
\node [arn_w] (5wc) at (11.2,0) {};

\begin{scope}[every node/.style={fill=white, circle, minimum width=0.01cm, minimum height=0.01cm}]
\path [->] (1) edge[fwd] (2);
\path [->] (2) edge[fwd] (3);
\path [->] (3) edge[fwd] (4);
\path [->] (4) edge[fwd] (5);
\path [->] (5w) edge[fwdds] (4w);
\path [->] (4w) edge[fwdds] (3w);
\path [->] (3w) edge[fwdds] (2w);
\path [->] (2w) edge[fwdds] (1w);
\path [->] (5wc) edge[Fwdds] (4wc);
\path [->] (4wc) edge[Fwdds] (3wc);
\path [->] (3wc) edge[Fwdds] (2wc);
\path [->] (2wc) edge[Fwdds] (1wc);

\end{scope}
\end{tikzpicture}
\label{fig:PartsSDDPalgrthms-APQP}
}

\subfigure[The adaptive partition-based SDDP algorithm with cautious cut generation (APCP-SDDP).]{
\begin{tikzpicture}
\node [arn_x] (1) at (0,1) {1};
\node [arn_x] (2) at (2.8,1) {2};
\node [arn_x] (3) at (5.6,1) {3};
\node [arn_x] (4) at (8.4,1) {4};
\node [arn_x] (5) at (11.2,1) {5};
\node [arn_w] (1w) at (0,0) {};
\node [arn_w] (2w) at (2.8,0) {};
\node [arn_w] (3w) at (5.6,0) {};
\node [arn_w] (4w) at (8.4,0) {};
\node [arn_w] (5w) at (11.2,0) {};
\node [arn_w] (1wc) at (0,-0.5) {};
\node [arn_w] (2wc) at (2.8,-0.5) {};
\node [arn_w] (3wc) at (5.6,-0.5) {};
\node [arn_w] (4wc) at (8.4,-0.5) {};
\node [arn_w] (5wc) at (11.2,-0.5) {};

\node [arn_w] (1w1) at (0,-1.05) {};
\node [arn_w] (2w1) at (2.8,-1.05) {};
\node [arn_w] (3w1) at (5.6,-1.05) {};
\node [arn_w] (4w1) at (8.4,-1.05) {};
\node [arn_w] (5w1) at (11.2,-1.05) {};
\node [arn_w] (1wc1) at (0,-0.505) {};
\node [arn_w] (2wc1) at (2.8,-0.505) {};
\node [arn_w] (3wc1) at (5.6,-0.505) {};
\node [arn_w] (4wc1) at (8.4,-0.505) {};
\node [arn_w] (5wc1) at (11.2,-0.505) {};

\begin{scope}[every node/.style={fill=white, circle, minimum width=0.01cm, minimum height=0.01cm}]
\path [->] (1) edge[fwd] (2);
\path [->] (2) edge[fwd] (3);
\path [->] (3) edge[fwd] (4);
\path [->] (4) edge[fwd] (5);
\path [->] (5w1) edge[Fwdds, bend left=50] (4w1);
\path [->] (4w1) edge[Fwdds, bend left=50] (3w1);
\path [->] (3w1) edge[Fwdds, bend left=50] (2w1);
\path [->] (2w1) edge[Fwdds, bend left=50] (1w1);
\path [->] (5wc1) edge[fwdds, bend left=50] (4wc1);
\path [->] (4wc1) edge[fwdds, bend left=50] (3wc1);
\path [->] (3wc1) edge[fwdds, bend left=50] (2wc1);
\path [->] (2wc1) edge[fwdds, bend left=50] (1wc1);

\path [->] (1w) edge[fwd, bend left=50] (2w);
\path [->] (2w) edge[fwd, bend left=50] (3w);
\path [->] (3w) edge[fwd, bend left=50] (4w);
\path [->] (4w) edge[fwd, bend left=50] (5w);
\path [->] (1wc) edge[fwd, bend left=50] (2wc);
\path [->] (2wc) edge[fwd, bend left=50] (3wc);
\path [->] (3wc) edge[fwd, bend left=50] (4wc);
\path [->] (4wc) edge[fwd, bend left=50] (5wc);

\end{scope}
\end{tikzpicture}
\label{fig:PartsSDDPalgrthms-APCP}
}

\subfigure[The adaptive partition-based SDDP with structured cut generation policies (SPAP-SDDP).]{
\begin{tikzpicture}
\node [arn_x] (1) at (0,1) {1};
\node [arn_r] (2) at (2.8,1) {2};
\node [arn_b] (3) at (5.6,1) {3};
\node [arn_b] (4) at (8.4,1) {4};
\node [arn_r] (5) at (11.2,1) {5};
\node [arn_w] (1w) at (0,0.5) {};
\node [arn_w] (2w) at (2.8,0.5) {};
\node [arn_w] (3w) at (5.6,0.5) {};
\node [arn_w] (4w) at (8.4,0.5) {};
\node [arn_w] (5w) at (11.2,0.5) {};
\node [arn_w] (1wc) at (0,0) {};
\node [arn_w] (2wc) at (2.8,0) {};
\node [arn_w] (3wc) at (5.6,0) {};
\node [arn_w] (4wc) at (8.4,0) {};
\node [arn_w] (5wc) at (11.2,0) {};

\begin{scope}[every node/.style={fill=white, circle, minimum width=0.01cm, minimum height=0.01cm}]
\path [->] (1) edge[fwd] (2);
\path [->] (2) edge[fwd] (3);
\path [->] (3) edge[fwd] (4);
\path [->] (4) edge[fwd] (5);
\path [->] (5w) edge[fwdds] (4w);
\path [->] (4w) edge[fwd] (3w);
\path [->] (3w) edge[fwd] (2w);
\path [->] (2w) edge[fwdds] (1w);
\path [->] (5wc) edge[Fwdds] (4wc);
\path [->] (2wc) edge[Fwdds] (1wc);

\end{scope}
\end{tikzpicture}
\label{fig:PartsSDDPalgrthms-SPAP}
}
\subfigure{
\begin{tikzpicture}
\matrix [draw,below left] at (current bounding box.south) {
\node [arn_xx,label=right: Stage $t$]at (-2,0) {$t$}; \\
\node [arn_rr,label=right:Stage of "less-importance": stage in a dry season] at (-2,0) {$t$}; \\
\node [arn_bb,label=right:Stage of "more-importance": stage in a wet season] at (-2,0) {$t$}; \\
\node [sqr_xx,label=right:Stage $t$ in the coarse tree $\hat{\T}(\N)$ ] at (-2,0) {}; \\
\node [arn_vsw] (0bd) at (-2,0) {};
\node [arn_vsw,label=right:Fine cut] (00bd) at (0.25,0) {};
\path [->] (00bd) edge[fwd, line width = 0.5 mm] (0bd);\\
\node [arn_vsw] (1bd) at (-2,0) {};
\node [arn_vsw,label=right: Forward pass] (11bd) at (0.25,0) {};
\path [->] (1bd) edge[fwd, line width = 0.5 mm] (11bd);\\
\node [arn_vsw] (3bd) at (-2,0) {};
\node [arn_vsw,label=right:Coarse cut] (33bd) at (0.25,0) {};
\path [->] (33bd) edge[fwdds,  line width = 0.5 mm, dash pattern = on 12pt] (3bd);\\
\node [arn_vsw] (4bd) at (-2,0) {};
\node [arn_vsw,label=right: Semi-coarse cut] (44bd) at (0.25,0) {};
\path [->] (44bd) edge[Fwdds,  line width = 0.5 mm] (4bd);\\
\node [arn_vsw] (6bd) at (-2,0) {};
\node [arn_vsw,label=right:Inner forward pass] (66bd) at (0.25,0) {};
\path [->] (6bd) edge[fwd, bend left=40,  line width = 0.5 mm] (66bd);\\
\node [arn_vsw] (7bd) at (-2,0) {};
\node [arn_vsw,label=right:Inner coarse cut] (77bd) at (0.25,0) {};
\path [->] (77bd) edge[fwdds, bend left=40, line width = 0.5 mm,  dash pattern = on 12pt] (7bd);\\
\node [arn_vsw,label=right:Inner semi-coarse cut] (77bd) at (0.25,0) {};
\path [->] (77bd) edge[Fwdds, bend left=40, line width = 0.5 mm] (7bd);\\
};
\end{tikzpicture}
}
\end{center}
\caption{Illustration of tree traversal strategies and various types of cutting planes used in the \textit{Refinement within SDDP} algorithms and the SPAP-SDDP.}
\label{fig:PartsSDDPalgrthms}
\end{figure}

Perhaps the most crucial question to this framework is how we can classify different stages into "more-important" and "less-important". The answer to this is clearly non-trivial, problem specific, and perhaps one that deserves a separate in-depth study by itself. Nonetheless, in our numerical experiments (see Section~\ref{sec:results}) we consider a heuristic approach of classification that we justify below.

We preface this by revisiting the hydro-thermal power generation planning problem that is considered in our numerical experiments. Readers are referred to, e.g., \cite{song2017level,de2017assessing}, for a more thorough discussion on the problem. In this problem, the decision maker aims to minimize the expected total cost which consists of: the power generation expenses and the penalty for the shortage in satisfying the demand. The stochasticity aspect of the problem arises due to the uncertainty about the amount of rainfall in the future -- which the decision maker can use to generate power via the interconnected network of hydro plants. As such, from a stochastic programming point of view, the amount of rain available at every stage $t=1, \dots, T$ is what defines as the valuable information that the decision maker will utilize in order to construct an optimal policy. In our implementation, we classify the stages as following: 
\begin{enumerate}
\item "Dry season" stage, which we label as "less important".
\item "Wet season" stage, which we label as "more important".
\end{enumerate}
The motivation for considering stages in the wet seasons to be of more importance and stages in the dry seasons to be of less importance is that, from an optimization point of view, the decision policy plays a more important role when it has more valuable resources at its disposal compared to when it does not. To put this into perspective, given the nature of the problem being a resource allocation/planning problem, an optimal policy is characterized by how balanced the availabilities of various types of resources are at times of abundance with their amounts at times of deficit. In dry seasons there is less of a decision to be made about resource allocation and more of a cost to be paid as a recourse action; whereas during the wet seasons the decision maker has to achieve a balance by allocating some of the available resources for generating power to meet the immediate demand while reserving some for hedging against the potential deficit in the future. 

In our implementation, we classify the dry and wet season stages by doing the following for every stage $t=2, \dots, T$:
\begin{enumerate} 
\item First, in order to differentiate between the different stages solely based on the inflow of rain, we equalize for everything by setting the water level carried over from the previous stage at every reservoir to be zero. That is, we assume that every reservoir is empty. This step is done by setting the state variable $x_{t-1} = 0$.  
\item Next, we optimize myopically with respect to stage $t$ by solving the stage $t$ subproblem for every realization of random vector $\xi_t \in \Xi_t$.
\item Let $z^{*}(t,\xi_t)$ be the optimal objective value corresponding to every respective problem solved in the previous step. 
\item Define $\bar{z}_t = \sum_{k=1}^{|\Xi_t|} z^{*}(t,\xi_{t,k})$ be the "worst-case" immediate cost at stage $t$, and let $\Z$ be a user prespecified parameter.
\item If $\bar{z}_t \leq \Z$, then stage $t$ is a "wet-stage". Otherwise, $t$ is a "dry-stage".
\end{enumerate}

\section{Numerical Results}
\label{sec:results}
In this section we report and analyze our numerical experiment results to show the empirical performances of the proposed algorithms. We first present the test instances and give an overview of different algorithms tested in the numerical experiments. Specifically, in Subsection~\ref{subsec:results_ROSDDP} we compare the \textit{Refinement outside SDDP} algorithms with the standard SDDP algorithm. In Subsection~\ref{subsec:results_RWSDDP} we compare the \textit{Refinement within SDDP} algorithms with the corresponding SDDP algorithms with different tree traversal strategies. Finally, in Subsection~\ref{subsec:results_SPSDDP} we compare the \textit{adaptive partition-based SDDP with structured cut generation policies} with the standard SDDP algorithm and different adaptive partition-based strategies.

We implemented all algorithms in \textit{Julia} 0.6.2, using package \textit{JuMP} 0.18.4 \cite{dunning2017jump}, with commercial solver \textit{Gurobi}, version 8.1.1 \cite{gurobi}. All the tests are conducted on Clemson University's primary high-performance computing cluster, the \textit{Palmetto} cluster, where we used an \textit{R830 Dell Intel Xeon} compute node with 2.60GHz and 1.0 TB memory. The number of cores is set to be 24.

\subsection{Test Instances and Algorithms}
\label{subsec:data_SDDPimplment}
As previously noted, we consider the multistage hydro-thermal power generation planning problem described in \cite{de2017assessing,song2017level}. We also use the same problem instance provided by E. Finardi and F. Beltr\'an, which models the Brazilian hydro-thermal power system. From the original data set, in order to create a variety of instances, we consider different planning horizons $T \in \{25, 61, 97, 120\}$ and sample sizes $|\Xi_t| \in \{50, 200, 1000\}$ of the random vector $\xi_t,\; \forall \; t=2, \dots T$. We let the number of realizations to be the same at every stage. For example, when $T=120$ and $|\Xi_t| =1000$, we will have $|\Xi_2| \times |\Xi_3| \times \dots |\Xi_{120}| = 1000^{119}$ scenarios (sample paths) .

Additionally, following on from the descriptions of  different tree traversal strategies in Subsection~\ref{subsubsec:tree_traversal} we implement two different variants of the SDDP algorithm:
\begin{enumerate}
\item \textit{SDDP with quick pass (SDDP-QP)}: adopt a quick pass strategy in traversing the scenario tree. We emphasize that SDDP-QP is the most commonly used variant of the \textit{SDDP} algorithm and we also refer to it as the \textit{standard SDDP} algorithm. \\
\item \textit{SDDP with cautious pass (SDDP-CP)}: adopt a cautious pass strategy in traversing the scenario tree.\\
\end{enumerate}
To get an initial lower bound for the cost-to-go function $\check \Q_t(\cdot), \; \forall \; t=1, 2, \dots T-1$, we solve the mean value problem with respect to the $(t+1)$-th stage problem by taking the expectation of the random vector $\xi_t$ and treating $\bar{x}_t$ in~\eqref{eq:1stagepbm} as decision variables. To measure the performances for different algorithms we report the lower bounds (LB) obtained by different algorithms after one, three and six hours (3600, 10800 and 21600 seconds, respectively) of processing. The LB progress is one of primary interests when solving an MSLP. We analyze the performance results by focusing on three factors:
\begin{enumerate}
\item The total number of stages in the planning horizon $T$.
\item The number of realizations per stage $|\Xi_t|$.
\item The processing time limit.
\end{enumerate}
In our experiments, we considered varying the number of sample paths at every iteration of the SDDP algorithm, however, we have found that using a single sample path per forward step (one scenario per iteration) works the best for all algorithms that we tested. 

\subsection{Numerical Results for the Refinement outside SDDP Strategy}
\label{subsec:results_ROSDDP}
\begin{figure}[htbp]
   \centering
  \subfigure[$T = 25, |\Xi_t|=50$. One hour of processing.]{
\includegraphics[scale=0.2175]{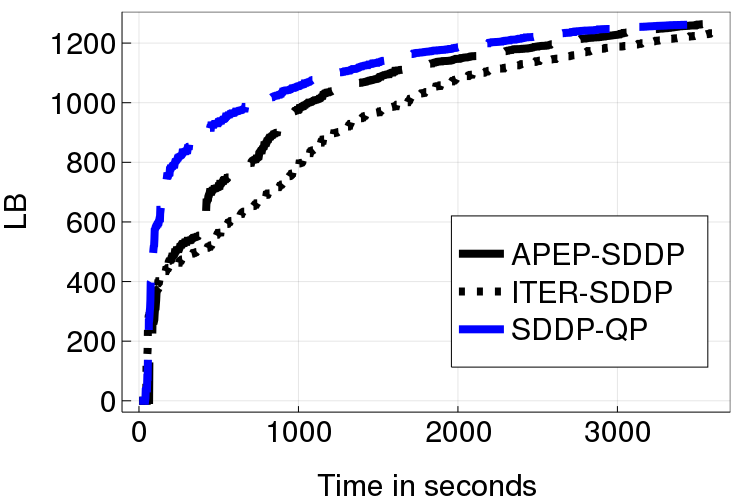}
  \label{fig:25T-50R-RO-1hrs} }
  \subfigure[$T = 25, |\Xi_t|=50$. Six hours of processing.]{
\includegraphics[scale=0.2175]{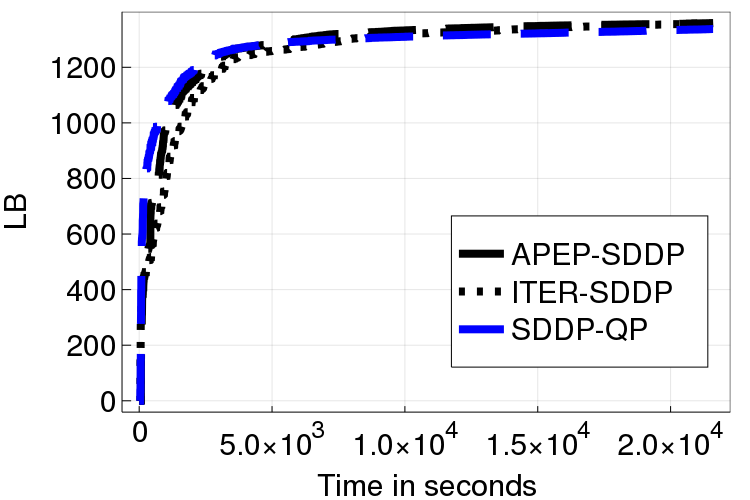}
  \label{fig:25T-50R-RO-6hrs} }
    \subfigure[$T = 61, |\Xi_t|=200$. One hour of processing.]{
\includegraphics[scale=0.2175]{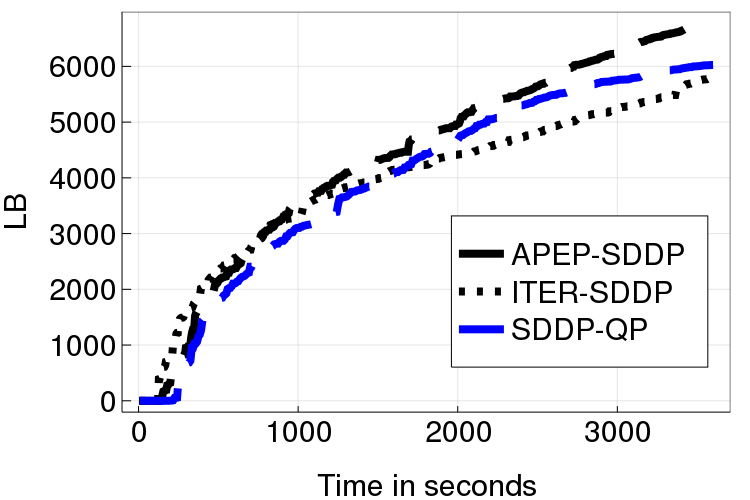}
  \label{fig:61T-200R-RO-1hrs} }
  \subfigure[$T = 61, |\Xi_t|=200$. Six hours of processing.]{
\includegraphics[scale=0.2175]{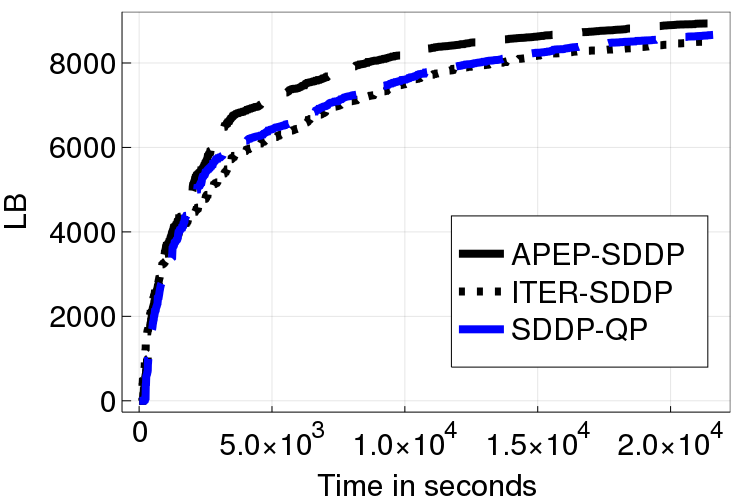}
  \label{fig:61T-200R-RO-6hrs} }
    \subfigure[$T = 97, |\Xi_t|=200$. One hour of processing.]{
\includegraphics[scale=0.2175]{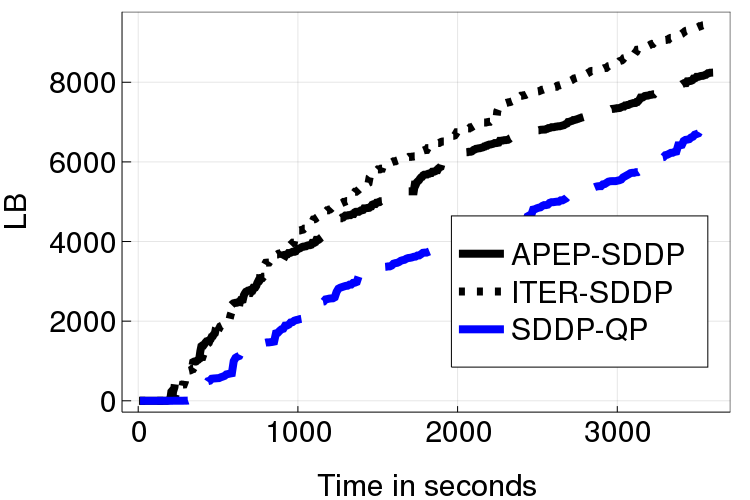}
  \label{fig:97T-200R-RO-1hrs} }
  \subfigure[$T = 97, |\Xi_t|=200$. Six hours of processing.]{
\includegraphics[scale=0.2175]{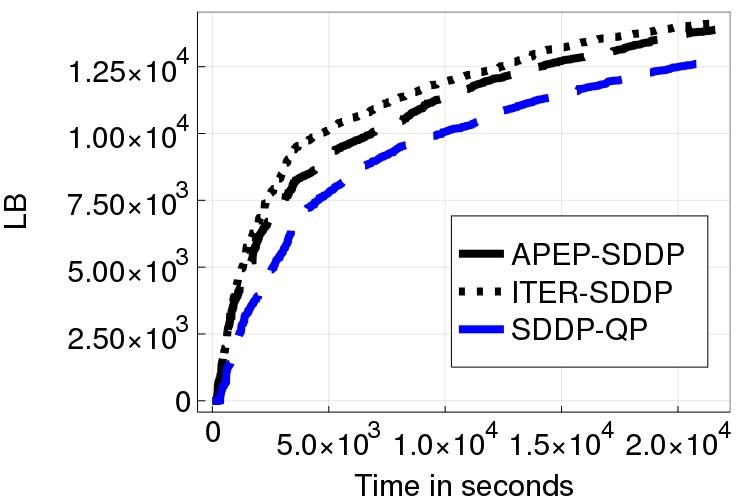}
  \label{fig:97T-200R-RO-6hrs} }
    \subfigure[$T = 120, |\Xi_t|=1000$. One hour of processing.]{
\includegraphics[scale=0.2175]{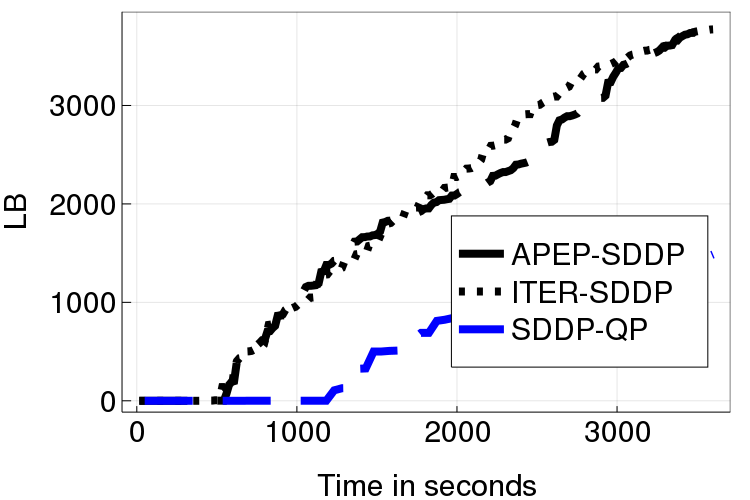}
  \label{fig:120T-1000R-RO-1hrs} }
  \subfigure[$T = 120, |\Xi_t|=1000$. Six hours of processing.]{
\includegraphics[scale=0.2175]{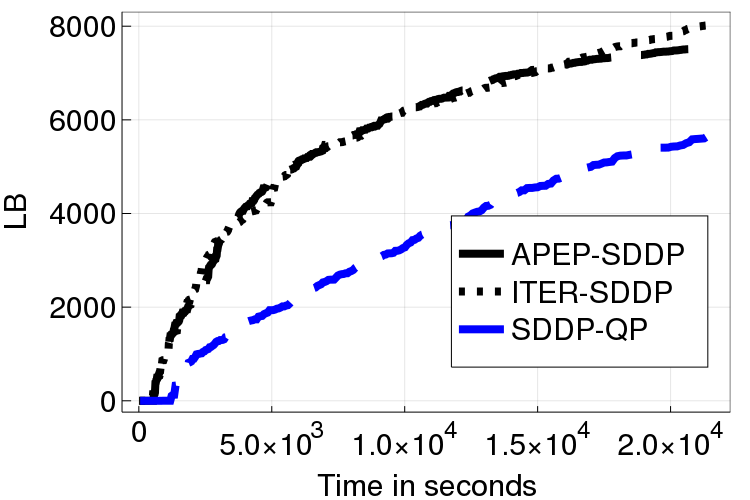}
  \label{fig:120T-1000R-RO-6hrs} }
  %
 \caption{Solution progress using \textit{Refinement outside SDDP} compared to SDDP-QP after one and six hours of processing on different instances.}
\label{fig:RO}
\end{figure}

We report and analyze the results of the two variants of the \textit{Refinement outside SDDP} strategy presented in Subsection~\ref{subsec:RefineOut}, namely APEP-SDDP and ITER-SDDP, and compare them to those obtained by the standard SDDP algorithm, i.e., SDDP-QP. We report the numerical results in Table~\ref{tab:RO-results} and a few selected instances in Figure~\ref{fig:RO}, from which we observe the following:
\begin{enumerate}
\item \underline{Performance with respect to $T$:}
\begin{itemize}
\item[$\bullet$] The overall performances of both APEP-SDDP and ITER-SDDP compared to that of the standard SDDP algorithm is consistently improving as the number of stages $T$ increases. Specifically, except for the case when $T=25$, both APEP-SDDP and ITER-SDDP outperform the SDDP-QP algorithm for all processing time limits.
\item[$\bullet$] In particular, this improvement in the performance is apparent when $T = 120$. For one, three and six hours of processing:
\begin{itemize}
\item APEP-SDDP outperformed SDDP-QP by $69\%, 32\%$ and $16\%$ on average.
\item ITER-SDDP outperformed SDDP-QP by $67\%, 36\%$ and $21\%$ on average.
\end{itemize}
\end{itemize}
\item \underline{Performance with respect to $|\Xi_t|$:}
\begin{itemize}
\item[$\bullet$] The overall performances of both APEP-SDDP and ITER-SDDP compared to that of the standard SDDP algorithm is consistently improving as $|\Xi_t|$ increases.
\item[$\bullet$] At its peak, when $|\Xi_t|= 1000$ and $T = 120$, the improvements over SDDP-QP reach to about $158\%$ and $154\%$ after one hour, $77\%$ and $77\%$ after three hours, and $34\%$ and $43\%$ after six hours, for APEP-SDDP and ITER-SDDP, respectively.
\end{itemize}
\item \underline{Performance with respect to the processing time limit:}
\begin{itemize}
\item[$\bullet$] Except for when $|\Xi_t| = 50$, the overall performances of the two algorithms APEP-SDDP and ITER-SDDP have the following trend under different processing time limits: as the processing time limit increases, the relative gap between the lower bounds of \textit{Refinement outside SDDP} algorithms and the SDDP algorithm, either deteriorates if it was superior or improves if it was inferior. 
\item[$\bullet$] For number of stages $T$ and realizations per stage $|\Xi_t|$ at the smaller end of our instances set, the \textit{Refinement outside SDDP} algorithms being inferior compared to the SDDP algorithm start to improve as processing time limit increases; whereas for larger instances, during earlier periods of processing time, \textit{Refinement outside SDDP} algorithms inherit a significant lead over the SDDP algorithm in its lower bound progress, but this wide margin steadily shrinks as the processing time increases.
\item[$\bullet$] In instances where the \textit{Refinement outside SDDP} algorithms are being outperformed by the SDDP algorithm, the advantage of the SDDP algorithm is at its highest after the first hour. 
\item[$\bullet$] The advantage of the \textit{Refinement outside SDDP} algorithms over the SDDP algorithm is most apparent after 3 hours of processing time limit, yielding a better performance in 7.5 out of 12 instances on average between the APEP-SDDP and ITER-SDDP algorithms.
\end{itemize}
\end{enumerate}
We attribute the aforementioned observations to the following:
\begin{itemize}
\item[$\bullet$] It should come as no surprise that, the larger the instance in terms of $T$ and $|\Xi_t|$, the better the performance of \textit{Refinement outside SDDP} algorithms should be. This advantage for the APEP-SDDP and ITER-SDDP algorithms over the standard SDDP algorithm, is a natural consequence due to the merits of adaptive partition-based strategies in making the cut generation effort adaptive to the solution progress.
\item[$\bullet$]  The aforementioned advantage does not hold for smaller instances, such as when $T=25$ and $|\Xi_t| =50$. This is because the computational savings provided by this framework via efficient cut generation effort from coarse scenario trees do not offset the significant inaccuracy inherited in these coarse cuts on these instances.    
This incompetence, however, steadily vanishes as the processing time limit increases when the coarse tree gets more and more refined.
\item[$\bullet$]  This explains the aforementioned observations regarding the performance with respect to the processing time limit:
\begin{itemize}
 \item In smaller instances, the larger the processing time limit is, the more the algorithm is able to compensate for the inaccuracy compromised during the early phase of the solution process.
 \item In larger instances, the larger the processing time limit is, the bigger that the size of the coarse tree gets, making the cut generation effort in the \textit{Refinement outside SDDP} algorithms similar to that of the standard SDDP algorithm. 
\end{itemize}
\item[$\bullet$] Finally, comparing between the two \textit{Refinement outside SDDP} algorithms, we see that their performances are comparable with $\textit{ITER-SDDP}$ prevailing in all instances except for when $T=25, |\Xi_t|=50$ and $T=120, |\Xi_t|=1000$. We attribute these two cases to the heuristically chosen parameters $n$ used in Algorithm~\ref{alg:APEP-SDDP}  as a criterion to perform partition refinement on ${\overrightarrow{\N}}^p$.
\end{itemize}
We conclude this subsection by emphasizing that, our work is an attempt to provide a framework which mitigates the computational burden of solving MSLPs brought by the \textit{curse-of-dimensionalty} due to the large number of stages $T$ and large number of realizations per stage $|\Xi_t|$. This, if anything, can only testify to the competitive nature of \textit{adaptive partition-based} strategies in solving large-scale problems. 
\begin{table}[]
\centering  	 																	
\begin{tabular}{|c|c|c|cc|cc|cc|c|}
\hline
\multicolumn{2}{|c|}{Instances}              & \multirow{2}{*}{Algorithm} & \multicolumn{2}{c|}{1 hour} & \multicolumn{2}{c|}{3hours} & \multicolumn{2}{c|}{6hours} & \multirow{2}{*}{\# iter} \\ \cline{1-2} \cline{4-9}
$T$                  & $|\Xi_t|$               &                            & LB             & \%LB        & LB            & \%LB        & LB            & \%LB        &                          \\ \hline
\multirow{6}{*}{25}  & \multirow{2}{*}{50}   & APEP  & 1265.2  & 0\%   & 1337.1  & 2\%   & 1358.9  & 2\%   & 4571 \\
                     &                       & ITER  & 1233.0  & -3\%  & 1322.6  & 1\%   & 1358.4  & 2\%   & 5428 \\ \cline{3-10}
                     & \multirow{2}{*}{200}  & APEP  & 552.3   & -5\%  & 625.6   & -4\%  & 671.2   & -2\%  & 2732 \\
                     &                       & ITER  & 573.6   & -1\%  & 669.5   & 3\%   & 713.1   & 4\%   & 3330 \\ \cline{3-10}
                     & \multirow{2}{*}{1000} & APEP  & 328.1   & -32\% & 454.0   & -17\% & 506.5   & -12\% & 1943 \\
                     &                       & ITER  & 371.7   & -23\% & 480.5   & -12\% & 529.2   & -8\%  & 2040 \\ \hline
\multirow{6}{*}{61}  & \multirow{2}{*}{50}   & APEP  & 11343.4 & 0\%   & 12848.7 & -1\%  & 13579.2 & -2\%  & 2506 \\
                     &                       & ITER  & 11636.3 & 2\%   & 12941.5 & 0\%   & 13640.4 & -1\%  & 2853 \\ \cline{3-10}
                     & \multirow{2}{*}{200}  & APEP  & 6775.8  & 12\%  & 8335.0  & 7\%   & 8943.7  & 3\%   & 1675 \\
                     &                       & ITER  & 5804.0  & -4\%  & 7686.7  & -1\%  & 8519.6  & -2\%  & 1921 \\ \cline{3-10}
                     & \multirow{2}{*}{1000} & APEP  & 3118.9  & 10\%  & 4578.3  & 1\%   & 5583.7  & 9\%   & 1070 \\
                     &                       & ITER  & 3427.5  & 21\%  & 5043.9  & 11\%  & 5858.3  & 15\%  & 1083 \\ \hline
\multirow{6}{*}{97}  & \multirow{2}{*}{50}   & APEP  & 19984.8 & -4\%  & 23342.6 & -6\%  & 25792.6 & -4\%  & 1993 \\
                     &                       & ITER  & 19998.9 & -4\%  & 23553.7 & -6\%  & 25697.9 & -4\%  & 2359 \\ \cline{3-10}
                     & \multirow{2}{*}{200}  & APEP  & 8242.9  & 19\%  & 11633.2 & 13\%  & 13882.0 & 9\%   & 1528 \\
                     &                       & ITER  & 9485.0  & 37\%  & 12136.4 & 18\%  & 14128.9 & 11\%  & 1475 \\ \cline{3-10}
                     & \multirow{2}{*}{1000} & APEP  & 3937.1  & 93\%  & 6628.2  & 37\%  & 8065.9  & 28\%  & 760  \\
                     &                       & ITER  & 3967.0  & 95\%  & 6001.5  & 24\%  & 8061.1  & 28\%  & 893  \\ \hline
\multirow{6}{*}{120} & \multirow{2}{*}{50}   & APEP  & 21890.5 & 0\%   & 27585.4 & -1\%  & 30639.1 & -3\%  & 1614 \\
                     &                       & ITER  & 22518.0 & 3\%   & 29393.6 & 6\%   & 32429.3 & 2\%   & 2101 \\ \cline{3-10}
                     & \multirow{2}{*}{200}  & APEP  & 9595.5  & 36\%  & 13956.9 & 20\%  & 16308.0 & 12\%  & 1339 \\
                     &                       & ITER  & 9578.7  & 36\%  & 13599.1 & 17\%  & 15664.3 & 7\%   & 1353 \\ \cline{3-10}
                     & \multirow{2}{*}{1000} & APEP  & 3836.4  & 158\% & 6345.8  & 77\%  & 7578.4  & 34\%  & 735  \\
                     &                       & ITER  & 3772.9  & 154\% & 6317.2  & 77\%  & 8065.5  & 43\%  & 729  \\ \hline
\end{tabular}
\caption{Lower bound progress obtained by the two different \textit{Refinement outside SDDP} algorithms~\ref{subsubsec:APEP-SDDP} and~\ref{subsubsec:ITER-SDDP} compared to the SDDP-QP algorithm. \%LB $= 100\times \frac{(LB-LB_{SDDP-QP})}{(LB_{SDDP-QP})}$.}\label{tab:RO-results}
\end{table}

\subsection{Numerical Results for the Refinement within SDDP Strategy}
\label{subsec:results_RWSDDP}
\begin{figure}[htbp]
   \centering
  \subfigure[$T = 25, |\Xi_t|=50$. One hour of processing.]{
\includegraphics[scale=0.2175]{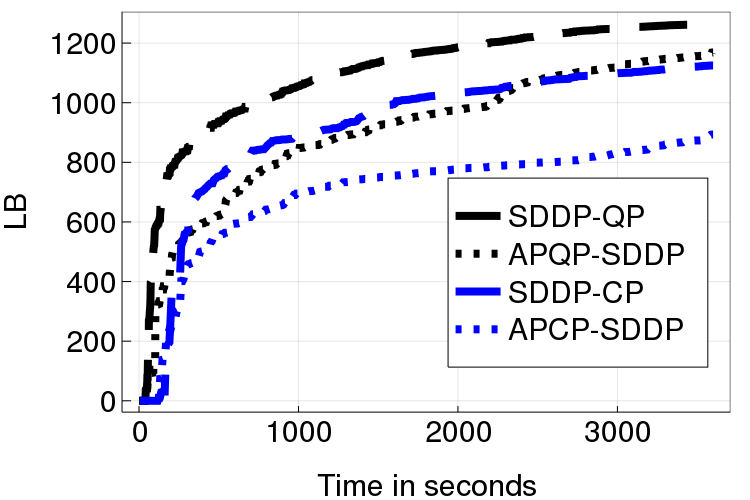}
  \label{fig:25T-50R-RW-1hrs} }
  \subfigure[$T = 25, |\Xi_t|=50$. Six hours of processing.]{
\includegraphics[scale=0.2175]{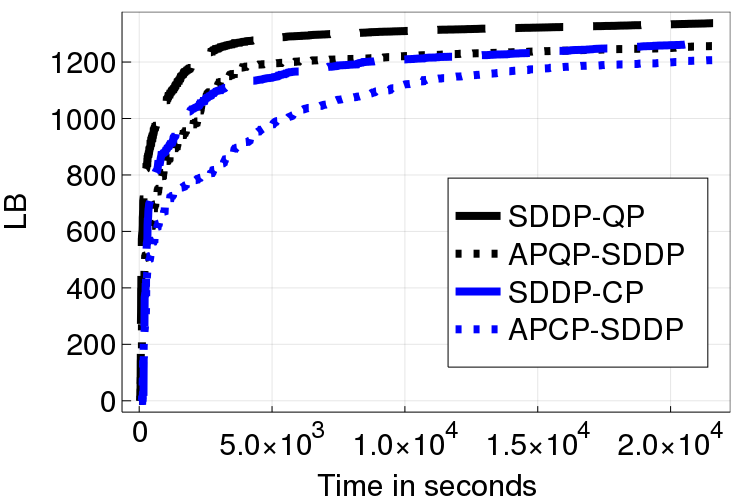}
  \label{fig:25T-50R-RW-6hrs} }
    \subfigure[$T = 61, |\Xi_t|=200$. One hour of processing.]{
\includegraphics[scale=0.2175]{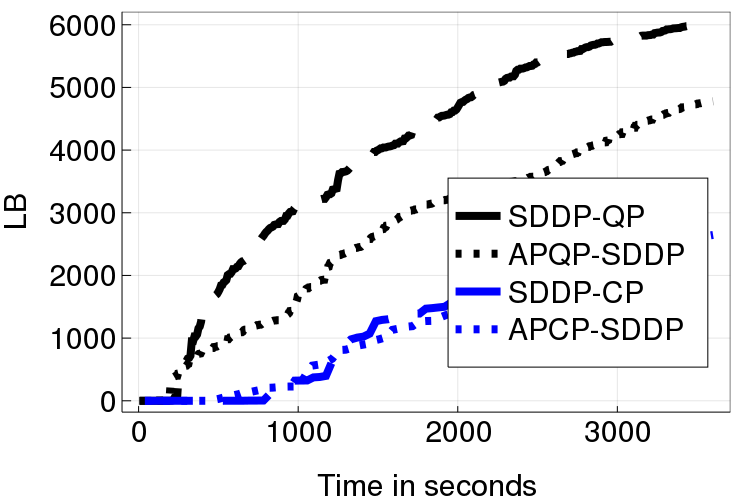}
  \label{fig:61T-200R-RW-1hrs} }
  \subfigure[$T = 61, |\Xi_t|=200$. Six hours of processing.]{
\includegraphics[scale=0.2175]{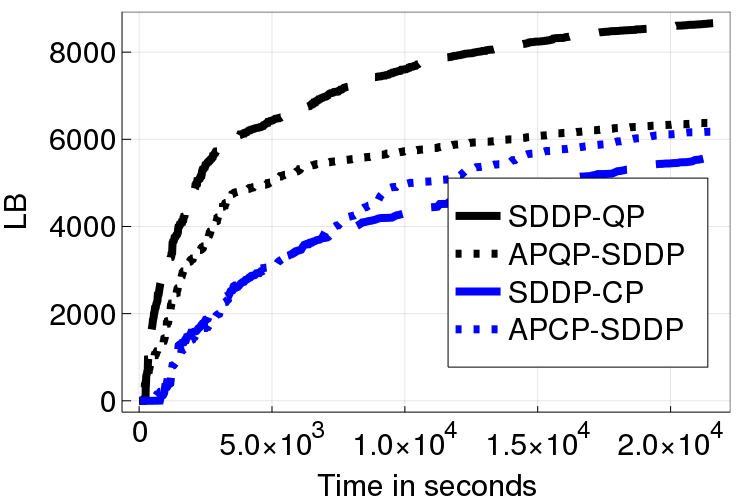}
  \label{fig:61T-200R-RW-6hrs} }
    \subfigure[$T = 97, |\Xi_t|=200$. One hour of processing.]{
\includegraphics[scale=0.2175]{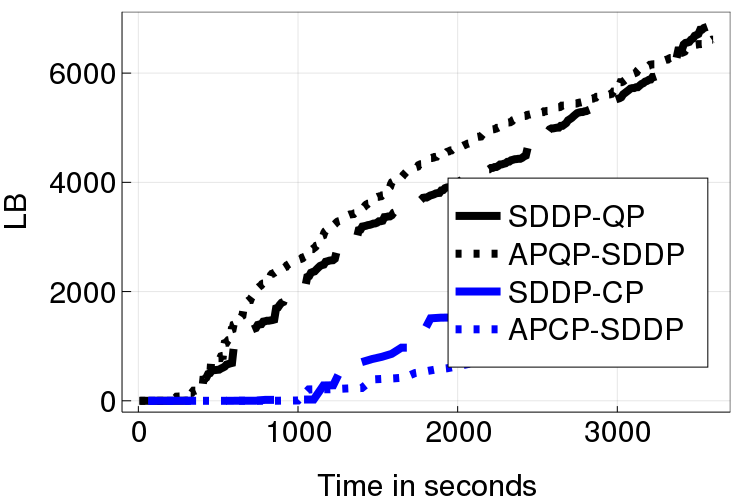}
  \label{fig:97T-200R-RW-1hrs} }
  \subfigure[$T = 97, |\Xi_t|=200$. Six hours of processing.]{
\includegraphics[scale=0.2175]{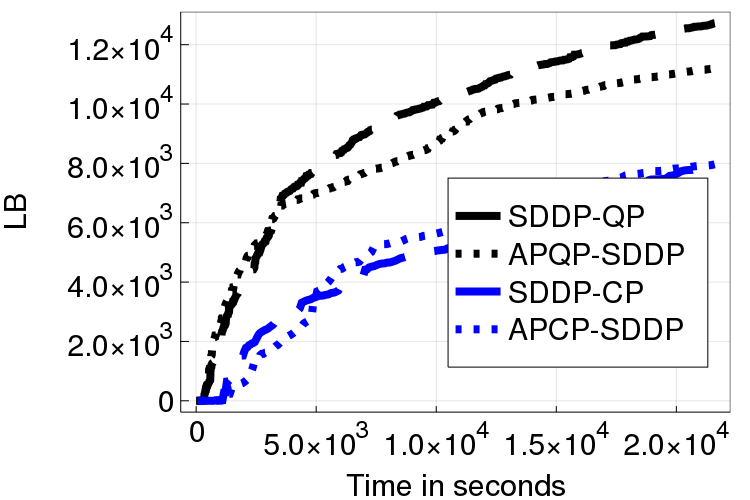}
  \label{fig:97T-200R-RW-6hrs} }
    \subfigure[$T = 120, |\Xi_t|=1000$. One hour of processing.]{
\includegraphics[scale=0.2175]{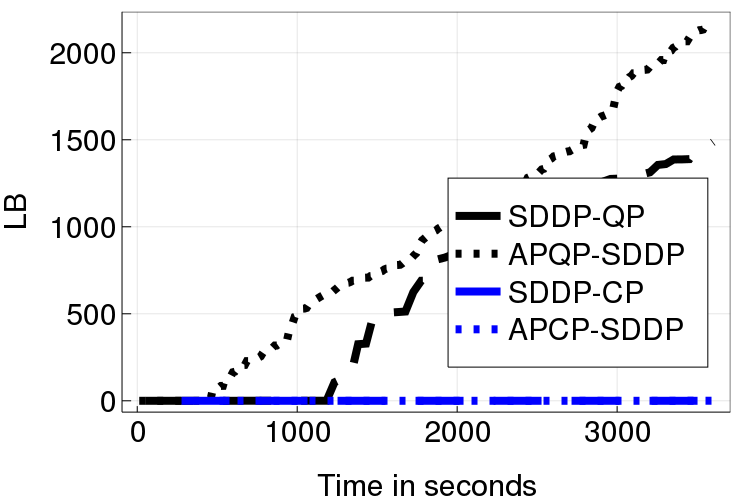}
  \label{fig:120T-1000R-RW-1hrs} }
  \subfigure[$T = 120, |\Xi_t|=1000$. Six hours of processing.]{
\includegraphics[scale=0.2175]{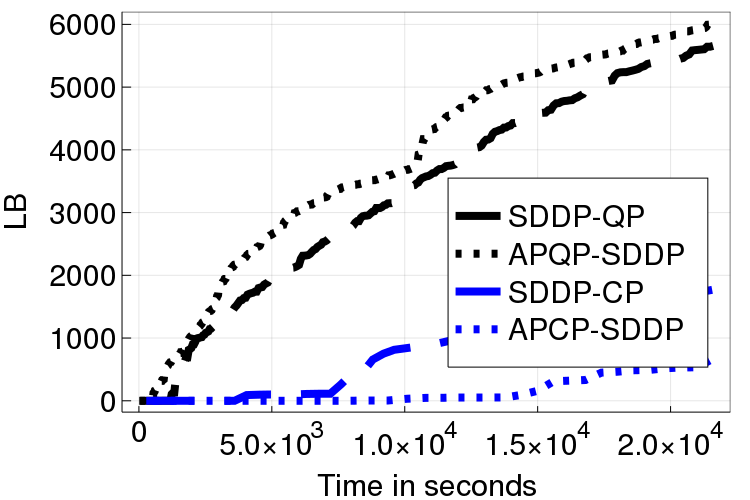}
  \label{fig:120T-1000R-RW-6hrs} }
  %
 \caption{Solution progress using \textit{Refinement within SDDP} algorithms compared to the corresponding SDDP algorithm under different tree traversal strategies for one and six hours of processing on different instances. \label{fig:RW}}
\end{figure}

In this subsection, we compare different \textit{Refinement within SDDP} algorithms described in Subsection~\ref{subsec:RefineWithin} and the corresponding SDDP algorithms (where no adaptive partition is employed) under different tree traversal strategies. We report the numerical results in Tables~\ref{tab:QP-results} and \ref{tab:CP-results} as well as a few selected instances in Figure~\ref{fig:RW}, from which we observe the following.
\begin{enumerate}
\item \underline{Performance with respect to $T$:}
\begin{itemize}
\item[$\bullet$] Similar to the observations made in Subsection~\ref{subsec:results_ROSDDP}, the performance of the APQP-SDDP algorithm compared to the SDDP-QP algorithm improves as $T$ increases. There is no consistent pattern comparing the performance of the APCP-SDDP algorithm with that of SDDP-CP as $T$ increases. 
\item[$\bullet$] Except for the performance of the APCP-SDDP algorithm in few instances of sizes in the mid range and the performance of APQP-SDDP in the large instances, the \textit{Refinement within SDDP} algorithms are consistently outperformed by the SDDP algorithm with the corresponding tree traversal strategy.  
\end{itemize}
\item \underline{Performance with respect to $|\Xi_t|$:}
\begin{itemize}
\item[$\bullet$] The overall results obtained by different algorithms do not have any obvious pattern associated with the number of realizations per stage $|\Xi_t|$.
\end{itemize}
\item \underline{Performance with respect to the processing time limit:}
\begin{itemize}
\item[$\bullet$] Similar to the observations made in Subsection~\ref{subsec:results_ROSDDP}, the overall performances of the two \textit{Refinement within SDDP}algorithms have the following trend under different processing time limits: as the processing time limit increases, the relative gap between the lower bounds obtained by the \textit{Refinement within SDDP} algorithms and the SDDP algorithm, either deteriorates if it was superior or improves if it was inferior. 
\end{itemize}
\end{enumerate}
We attribute the aforementioned observations to the following:
\begin{itemize}
\item[$\bullet$] As one might expect, allocating a significant cut generation effort to a small number of sample paths, is surely a waste of computational budget. This perhaps justifies the overall disappointing results of any algorithm with a tree traversal strategy that has some cautiousness aspects to it, since a cautious tree traversal strategy usually results in a decision policy that is overfitting to the subset of sample path(s) visited so far in the solution process. 
\item[$\bullet$] This overfitting in the decision policy could also explain the absence for any clear pattern in the performance of APCP-SDDP compared to that of SDDP-CP in different instances, since the performances depend on the sample paths visited by the algorithms, which are randomly generated.
\item[$\bullet$] When analyzing the incompetence in the performance of \textit{Refinement within SDDP} algorithms, and in particular the APQP-SDDP algorithm compared to \textit{Refinement outside SDDP} algorithms, it is important to note the following:
\begin{itemize}
\item In the \textit{Refinement outside SDDP} algorithms, where we only generate coarse cuts by implementing the standard SDDP algorithm on the coarse tree $\hat{\T}({\overrightarrow{\N}}^p)$, there is a criterion by which we measure the added value of these coarse cuts to the current cost-to-go function approximations $\check \Q^{p}_t(\cdot)$. This is done by keeping track of the lower bound progress when generating cuts from the coarse tree at every iteration (see Algorithm~\ref{alg:APEP-SDDP}). Hence, we restrict the algorithm to allocate some of the computational budget to generating coarse cuts using ${\overrightarrow{\N}}^p$, only if these coarse cuts have significant added value to them. Otherwise, ${\overrightarrow{\N}}^p$ will be refined. 
\item In the \textit{Refinement within SDDP} algorithms, the process is continuously attempting to generate coarse cuts, irrespective of their added value to the decision policy. This might hinder the performance of the algorithm, especially in the case of APQP-SDDP algorithm, where the existence of any violated cut, regardless of its quality, is the criterion for going back from stage $t$ to $t-1$. This makes refining the coarse tree (in order to obtain some accuracy) less frequent. 
\end{itemize}
\end{itemize}
We conclude this subsection by mentioning that, while the coarse cuts generated by the \textit{Refinement within SDDP} algorithms do not seem to serve its desired purpose, as far as the previous analysis goes, this may not necessarily be the case under different time limits or different integration scheme such as the SPAP-SDDP algorithm that we shall discuss next. 
\begin{table}[htbp]
\centering  																			
\caption{Lower bound progress obtained by the APQP-SDDP~\eqref{subsubsec:APQP-SDDP} algorithm compared to the SDDP-QP algorithm. \%LB $= 100\times \left(\frac{LB_{APQP-SDDP}-LB_{SDDP-QP}}{LB_{SDDP-QP}}\right)$.}
 \label{tab:QP-results}
\begin{tabular}{|c|c|cc|cc|cc|c|}
\hline
\multicolumn{2}{|c|}{Instances} & \multicolumn{2}{c|}{1 hour} & \multicolumn{2}{c|}{3hours} & \multicolumn{2}{c|}{6hours} & \multirow{2}{*}{\# iter} \\ \cline{1-8}
$T$                   & $|\Xi_t|$ & LB             & \%LB        & LB            & \%LB        & LB            & \%LB        &                          \\ \hline
\multirow{3}{*}{25}   & 50      & 1168.6         & -8\%        & 1224.5        & -7\%        & 1256.6        & -6\%        & 5697                     \\ 
                      & 200     & 439.0          & -25\%       & 534.0         & -18\%       & 587.9         & -14\%       & 3480                     \\  
                      & 1000    & 172.3          & -64\%       & 248.7         & -54\%       & 287.9         & -50\%       & 1875                     \\ \hline
\multirow{3}{*}{61}   & 50      & 9796.3         & -14\%       & 11240.7       & -13\%       & 11956.9       & -13\%       & 3417                     \\ 
                      & 200     & 4778.4         & -21\%       & 5771.7        & -26\%       & 6384.7        & -26\%       & 2129                     \\  
                      & 1000    & 2027.5         & -28\%       & 3056.9        & -33\%       & 3745.6        & -27\%       & 1254                     \\ \hline
\multirow{3}{*}{97}   & 50      & 15549.1        & -26\%       & 18132.6       & -27\%       & 19809.2       & -26\%       & 2647                     \\ 
                      & 200     & 6613.2         & -4\%        & 9137.5        & -11\%       & 11197.4       & -12\%       & 1628                     \\ 
                      & 1000    & 2058.1         & 1\%         & 5139.5        & 6\%         & 6327.5        & 1\%         & 854                      \\ \hline
\multirow{3}{*}{120}  & 50      & 17859.8        & -18\%       & 23961.7       & -14\%       & 27449.8       & -13\%       & 2392                     \\  
                      & 200     & 5886.7         & -17\%       & 9040.6        & -22\%       & 11962.0       & -18\%       & 1558                     \\
                      & 1000    & 2170.2         & 46\%        & 4204.8        & 18\%        & 6019.7        & 6\%         & 810                      \\ \hline
\end{tabular}
\end{table}

\begin{table}[htbp]
\centering 																					    
\caption{Lower bound progress obtained by the APCP-SDDP~\eqref{subsubsec:APCP-SDDP} algorithm compared to the SDDP-CP algorithm. \%LB $= 100\times \left(\frac{LB_{APCP-SDDP}-LB_{SDDP-CP}}{LB_{SDDP-CP}}\right)$.}
 \label{tab:CP-results}
\begin{tabular}{|c|c|cc|cc|cc|c|}
\hline
\multicolumn{2}{|c|}{Instances} & \multicolumn{2}{c|}{1 hour} & \multicolumn{2}{c|}{3hours} & \multicolumn{2}{c|}{6hours} & \multirow{2}{*}{\# iter} \\ \cline{1-8}
$T$                   & $|\Xi_t|$ & LB            & \%LB         & LB            & \%LB        & LB            & \%LB        &                          \\ \hline
\multirow{3}{*}{25}   & 50      & 893.2         & -21\%        & 1136.7        & -6\%        & 1206.5        & -5\%        & 1282                     \\ 
                      & 200     & 375.6         & -22\%        & 514.0         & -8\%        & 578.4         & -9\%        & 628                      \\ 
                      & 1000    & 174.3         & -36\%        & 351.6         & -5\%        & 392.2         & -4\%        & 172                      \\ \hline
\multirow{3}{*}{61}   & 50      & 8815.5        & 4\%          & 10851.3       & -4\%        & 12118.3       & -2\%        & 671                      \\  
                      & 200     & 2646.0        & 2\%          & 5028.6        & 14\%        & 6173.9        & 10\%        & 303                      \\ 
                      & 1000    & 0.5           & -100\%       & 1520.9        & -28\%       & 2962.7        & -7\%        & 85                       \\ \hline
\multirow{3}{*}{97}   & 50      & 13227.7       & 18\%         & 18933.9       & 8\%         & 21521.8       & 4\%         & 498                      \\ 
                      & 200     & 2004.0        & -29\%        & 5954.4        & 15\%        & 7969.6        & 0\%         & 215                      \\ 
                      & 1000    & 0.0           & -100\%       & 219.9         & -78\%       & 1258.9        & -58\%       & 44                       \\ \hline
\multirow{3}{*}{120}  & 50      & 11195.0       & -7\%         & 19532.9       & 0\%         & 23950.4       & 1\%         & 445                      \\  
                      & 200     & 1443.6        & -38\%        & 4752.1        & -10\%       & 6604.3        & -12\%       & 175                      \\  
                      & 1000    & 0.2           & 3304\%       & 46.5          & -95\%       & 672.6         & -62\%       & 33                       \\ \hline
\end{tabular}
\end{table}

\subsection{Computational Experiments on the Adaptive Partition-based SDDP with Structured Cut Generation Policies}
\label{subsec:results_SPSDDP}
\begin{figure}[htbp]
   \centering
  \subfigure[$T = 25, |\Xi_t|=50$. One hour of processing.]{
\includegraphics[scale=0.2175]{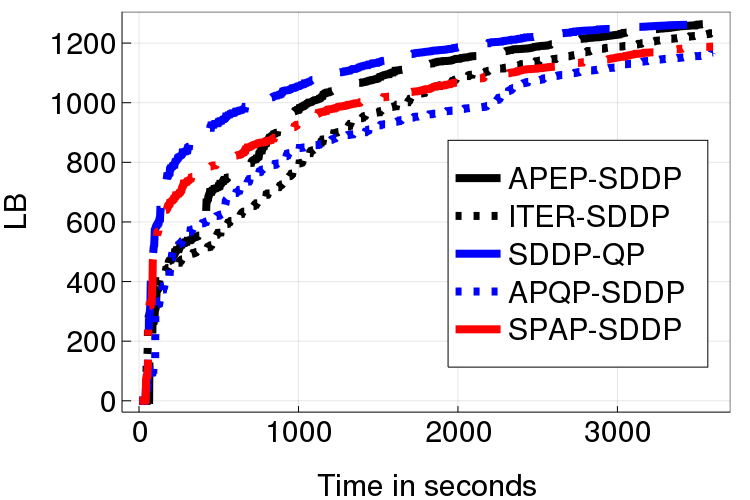}
  \label{fig:25T-50R-SPAP-1hrs} }
  \subfigure[$T = 25, |\Xi_t|=50$. Six hours of processing.]{
\includegraphics[scale=0.2175]{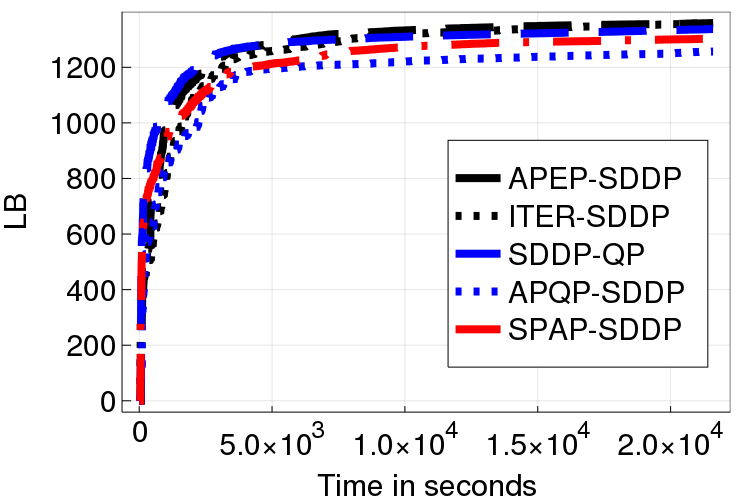}
  \label{fig:25T-50R-SPAP-6hrs} }
    \subfigure[$T = 61, |\Xi_t|=200$. One hour of processing.]{
\includegraphics[scale=0.2175]{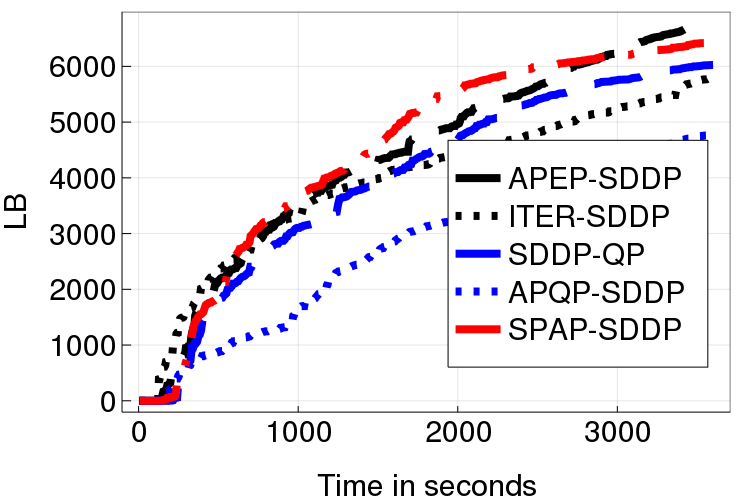}
  \label{fig:61T-200R-SPAP-1hrs} }
  \subfigure[$T = 61, |\Xi_t|=200$. Six hours of processing.]{
\includegraphics[scale=0.2175]{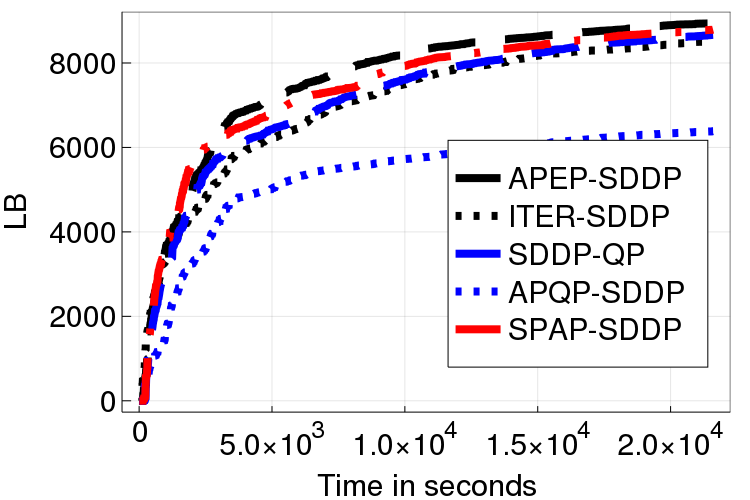}
  \label{fig:61T-200R-SPAP-6hrs} }
    \subfigure[$T = 97, |\Xi_t|=200$. One hour of processing.]{
\includegraphics[scale=0.2175]{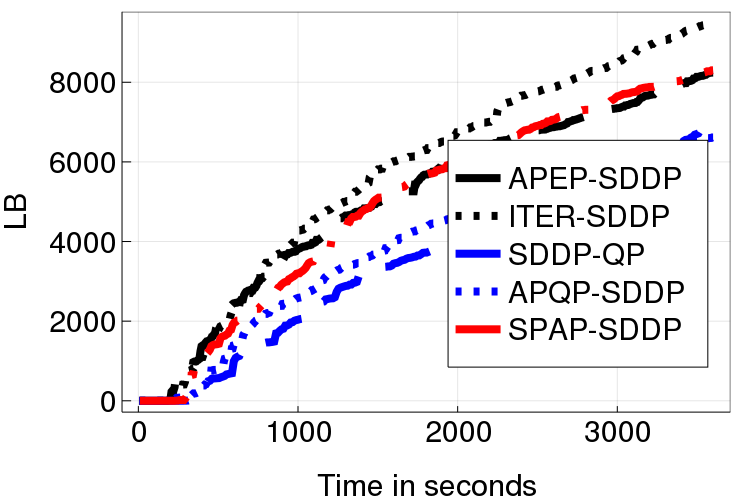}
  \label{fig:97T-200R-SPAP-1hrs} }
  \subfigure[$T = 97, |\Xi_t|=200$. Six hours of processing.]{
\includegraphics[scale=0.2175]{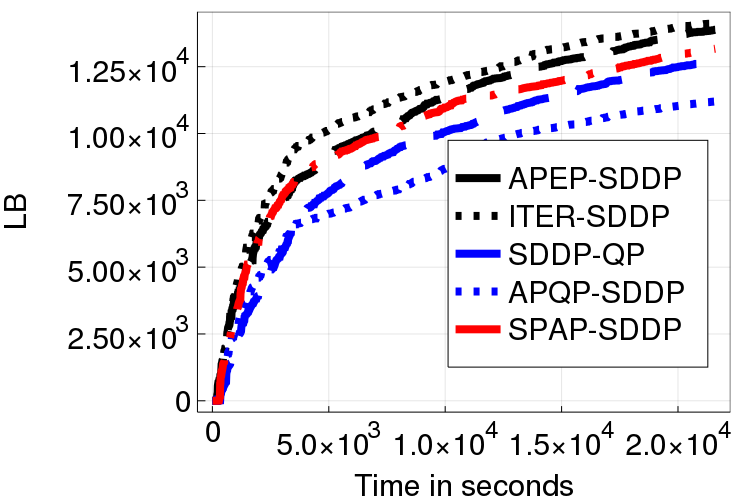}
  \label{fig:97T-200R-SPAP-6hrs} }
    \subfigure[$T = 120, |\Xi_t|=1000$. One hour of processing.]{
\includegraphics[scale=0.2175]{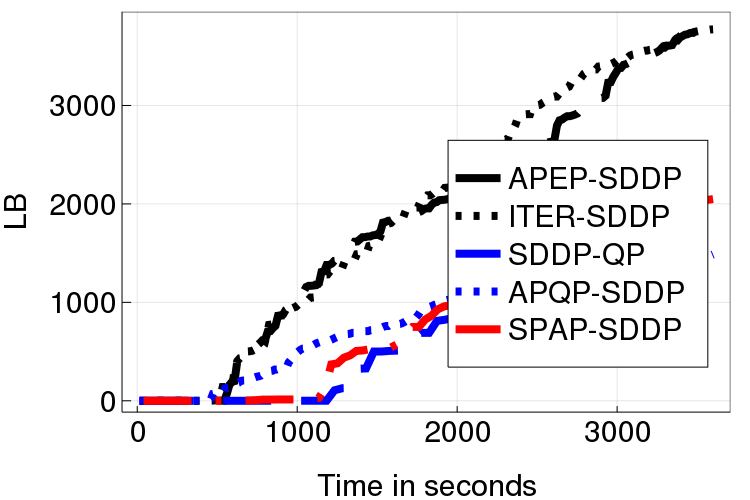}
  \label{fig:120T-1000R-SPAP-1hrs} }
  \subfigure[$T = 120, |\Xi_t|=1000$. Six hours of processing.]{
\includegraphics[scale=0.2175]{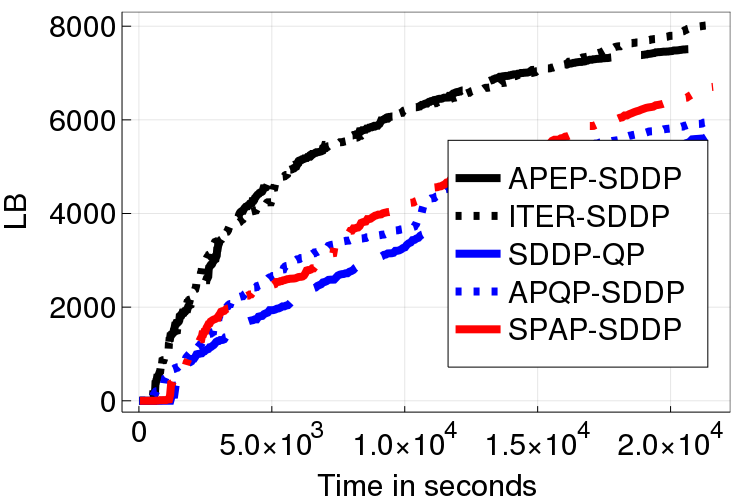}
  \label{fig:120T-1000R-SPAP-6hrs} }
  %
  \caption{Solution progress using SPAP-SDDP algorithm compared to \textit{Refinement outside SDDP} algorithms, APQP-SDDP and SDDP-QP after one and the six hours of processing on different instances.}
\label{fig:SPAP}
\end{figure}

In this subsection we report and analyze the results for the performance of the SPAP-SDDP algorithm presented in Subsection~\ref{subsubsec:SPAP-SDDP}. In Table~\ref{tab:SPAP-SDDP-results} we report its performance and compare it to the standard SDDP algorithm. Additionally, we illustrate a few selected instances in Figure~\ref{fig:SPAP} which compare the performance of the SPAP-SDDP algorithm with that of different \textit{Refinement outside SDDP} algorithms, as well as the APQP-SDDP and SDDP-QP algorithms.

The overall pattern in the behavior of the SPAP-SDDP algorithm for different number of stages $T$, realizations per stage $|\Xi_t|$ and processing time limits, is very similar to the observations made regarding the \textit{Refinement outside SDDP} algorithms. We summarize these observations as follows:
\begin{itemize}
\item[$\bullet$] Except for when $T=25$, the performance of the SPAP-SDDP algorithm compared to that of the standard SDDP algorithm is consistently improving as $T$ and/or $|\Xi_t|$ increases.
\item[$\bullet$] Similar to the observations made in Subsection~\ref{subsec:results_ROSDDP}, the overall performances of the SPAP-SDDP algorithm have the following trend under different processing time limits: as the processing time limit increases, the relative gap between the lower bounds of the SPAP-SDDP algorithm and the SDDP algorithm, either deteriorates if it was superior or improves if it was inferior. 

\end{itemize}
We attribute the aforementioned observations to the following:
\begin{itemize}
\item[$\bullet$] Most of the reasoning made in Subsection~\ref{subsec:results_ROSDDP} regarding the performance of the \textit{Refinement outside SDDP} algorithms can also be made to the SPAP-SDDP algorithm. This reasoning being, in large instances, where accuracy is more computationally expensive to obtain, the computational savings come from the fact that the algorithm makes the cut generation effort adaptive to the solution progress. 
Except here, this adaptability is not integrated by the value which a coarse cut adds immediately to the lower bound progress, but instead, by the added value of a coarse cut from a particular stage to the decision policy at giving point in the processing time, which affects the lower bound progress implicitly.  
\item[$\bullet$] The most notable difference between the performance the SPAP-SDDP algorithm and the \textit{Refinement outside SDDP} algorithms is that, the performance of the SPAP-SDDP algorithm is very stable in outperforming SDDP-QP compared to that of the \textit{Refinement outside SDDP} algorithms. However, when \textit{Refinement outside SDDP} algorithms do outperform SDDP-QP, they outperform by a large margin. 
\item[$\bullet$] Overall, it is not difficult to see that the SPAP-SDDP algorithm bridges the gap between the computational ease of generating excess of inaccurate coarse cuts using \textit{APQP-SDDP} and the computational burden of generating a few, but accurate fine cuts using the standard SDDP algorithm. 
\end{itemize}

\begin{table}[htbp]
\centering 																			
\caption {Lower bound progress obtained by the SPAP-SDDP~\eqref{subsubsec:SPAP-SDDP} algorithm compared to the SDDP-QP algorithm. \%LB $= 100\times \left(\frac{LB_{SPAP-SDDP}-LB_{SDDP-QP}}{LB_{SDDP-QP}}\right)$.}
 \label{tab:SPAP-SDDP-results}
\begin{tabular}{|c|c|cc|cc|cc|c|}
\hline
\multicolumn{2}{|c|}{Instances} & \multicolumn{2}{c|}{1 hour} & \multicolumn{2}{c|}{3hours} & \multicolumn{2}{c|}{6hours} & \multirow{2}{*}{\# iter} \\ \cline{1-8}
$T$                   & $|\Xi_t|$ & LB             & \%LB        & LB            & \%LB        & LB             & \%LB       &        \\ \hline
\multirow{3}{*}{25}   & 50      & 1188.4         & -6\%        & 1277.3        & -3\%        & 1302.7         & -3\%       & 4908                   \\ 
                      & 200     & 556.4          & -4\%        & 643.6         & -1\%        & 703.1          & 3\%        & 2671   \\ 
                      & 1000    & 436.1          & -10\%       & 493.6         & -10\%       & 533.0          & -8\%       & 1084 \\ \hline
\multirow{3}{*}{61}   & 50      & 11950.2        & 5\%         & 13030.3       & 0\%         & 13707.9        & -1\%       & 2704 \\ 
                      & 200     & 6434.7         & 7\%         & 8062.1        & 4\%         & 8790.2         & 1\%        & 1448 \\ 
                      & 1000    & 2716.1         & -4\%        & 4539.8        & 0\%         & 5344.5         & 5\%        & 561 \\ \hline
\multirow{3}{*}{97}   & 50      & 21178.3        & 1\%         & 25112.7       & 1\%         & 27129.1        & 1\%        & 2094 \\ 
                      & 200     & 8301.5         & 20\%        & 11283.8       & 10\%        & 13172.5        & 4\%        & 1088 \\
                      & 1000    & 3355.3         & 65\%        & 5944.7        & 23\%        & 7420.3         & 18\%       & 415 \\ \hline
\multirow{3}{*}{120}  & 50      & 23801.5        & 9\%         & 29635.3       & 7\%         & 32543.9        & 3\%        & 1897 \\ 
                      & 200     & 8391.5         & 19\%        & 12564.2       & 8\%         & 15171.2        & 4\%        & 978 \\
                      & 1000    & 2049.4         & 38\%        & 4449.3        & 24\%        & 6704.2         & 19\%       & 358 \\ \hline
\end{tabular}
\end{table}

\paragraph{Summary of numerical experiment results.} In sum, based on the aforementioned analysis on the numerical results, the \textit{Refinement outside SDDP} algorithms and the SPAP-SDDP algorithm yield the most significant improvements over the standard SDDP algorithm in terms of the lower bound progress. We thereby suggest that as a rule of thumb, integrating the adaptive partition-based strategies into the SDDP algorithm should be done via the \textit{Refinement outside SDDP} approach, and structured cut generation policies should be pursued (like the SPAP-SDDP algorithm) if any structure can be exploited from the underlying problem instance.

\section{Conclusions}
\label{sec:conclusions}
In this study, we have investigated various ways to enhance the performance of the SDDP algorithm in terms of its lower bound progress by employing various inexact cut generations and scenario-tree traversal strategies. Specifically, we have integrated the adaptive partition-based approaches, which have been shown to be effective in two-stage stochastic programs, into the SDDP algorithm for multi-stage stochastic programs in two different manners: performing partition refinement within the SDDP and outside the SDDP algorithm. In addition, we have proposed a structured cut generation strategy across all stages, which takes advantage of the underlying seasonal uncertainty structure in the class of problems that we use as the test instances. We have conducted extensive numerical experiments to empirically validate the effectiveness of the proposed algorithms and compare them to the standard SDDP algorithm. We have identified several directions to pursue for future research. First, it is of interest to investigate how the proposed algorithms can be applied to address more challenging problem classes such as distributionally robust multistage stochastic programs and multistage stochastic integer programs. Second, from an algorithmic perspective it is worth investigating novel cut generation strategies that are adaptive to the underlying problem structure such as the decision policy structures and/or uncertainty structures during the solution process, as opposed to imposing these structures a priori as we did in the SPAP-SDDP algorithm.

\section{acknowledgements}
\label{acknw}
Clemson University is acknowledged for generous allotment of compute time on Palmetto cluster. The authors acknowledge partial support by the National Science Foundation [Grant CMMI 1854960]. Any opinions, findings, and conclusions or recommendations expressed in this material are those of the authors and do not necessarily reflect the views of the National Science Foundation.
\bibliographystyle{unsrt}  
\bibliography{references}   



\end{document}